\documentclass[12pt,a4paper]{amsart}
\usepackage[dvipsnames]{xcolor}
\usepackage{amsthm,amsmath}
\usepackage{graphicx}
\usepackage{caption}
\usepackage{subcaption}
\usepackage{lipsum}
\usepackage{mathtools}
\usepackage{booktabs}
\mathtoolsset{showonlyrefs}

\graphicspath{ {./figures/} }

\theoremstyle{definition}
\newtheorem{definition}{Definition}[section]

\theoremstyle{plain}

\def \dd {\mathrm{d}}
\newcommand{\e}{\mathrm{e}}

\usepackage{geometry}
\usepackage{mathtools}
\usepackage{leftindex}
\usepackage{comment}
\usepackage{bm}

\usepackage{lineno}
\usepackage{mathtools}
\mathtoolsset{showonlyrefs}

\definecolor{cornellred}{rgb}{0.7, 0.11, 0.11}

\usepackage[colorlinks=true,urlcolor=RoyalBlue,citecolor=cornellred,linkcolor=RoyalBlue,linktocpage,pdfpagelabels,
bookmarksnumbered,bookmarksopen]{hyperref}

\title[Exploring Exponential Runge-Kutta Methods]{Exploring Exponential Runge-Kutta Methods: \\A Survey}

\author[A.\,And\`o ]{Alessia And\`o}
\author[N.\,Cangiotti]{Nicol\`o Cangiotti}
\author[M.\,Sensi]{Mattia Sensi}

\address[A.\,And\`o]{Department of Mathematics, Computer Science and Physics\newline\indent University of Udine \newline\indent
via delle scienze 206, 33100 Udine, Italy}
\email{alessia.ando@uniud.it}

\address[N.\,Cangiotti]{Department of Mathematics\newline\indent Politecnico di Milano \newline\indent
via Bonardi 9, Campus Leonardo, 20133 Milan, Italy}
\email{nicolo.cangiotti@polimi.it}

\address[M.\,Sensi]{Department of Mathematics \newline\indent University of Trento\newline\indent
Via Sommarive 14, 38123 Povo (Trento), Italy}
\email{mattia.sensi@unitn.it}


\begin{document}
\begin{abstract} 
In this survey, we provide an in-depth investigation of exponential Runge-Kutta methods for the numerical integration of initial-value problems. These methods offer a valuable synthesis between classical Runge-Kutta methods, introduced more than a century ago, and exponential integrators, which date back to the 1960s. This manuscript presents both a historical analysis of the development of these methods up to the present day and several examples aimed at making the topic accessible to a broad audience. 
\end{abstract}

\maketitle

\section{Introduction}\label{Sec1}
Numerical integration of systems of Ordinary Differential Equations (ODEs) is one of the fundamental pillars of numerical analysis. In particular, many models stemming from real-world scenarios, with applications ranging from quantum mechanics to fluid dynamics, naturally lead to systems of the form: 
\begin{equation}
\label{eq:start}
u'(t)=Au(t)+g(t,u(t)), \quad u(0)=u_0,  
\end{equation}
where $u(t)\in\mathbb{R}^d$, $t\geq 0$, $A$ is a linear operator (typically a matrix) describing the (typically stiff)  part of the problem, while $g$ represents the (typically non-stiff) nonlinear component. Remarkably, despite its fundamental importance, there is no univocal and universally accepted definition of stiffness. Some mathematicians characterize it through the presence of components of the solution that evolve on very different time scales, others define it in terms of the ratio between the eigenvalues of the Jacobian matrix of the system, while still others relate it to the practical difficulties that arise when using explicit numerical methods. This multiplicity of definitions reflects the intrinsic complexity of the phenomenon and the different perspectives from which it can be analyzed. We shall provide a more exhaustive discussion about the topic in Sect. \ref{Sec21}. 

Runge-Kutta (RK) methods represent a milestone in the development of numerical methods for ordinary differential equations. Initially introduced by Carl Runge in 1895 \cite{runge1895}
and later refined by Martin Wilhelm Kutta in 1901 \cite{kutta1901}, these methods revolutionized the approach to the numerical solution of ODEs. Their classical formulation, based on a linear combination of function evaluations at appropriately chosen intermediate points, offers a remarkable balance between accuracy and simplicity of implementation. However, despite their historical success, classical \emph{explicit} RK methods have shown significant limitations when applied to more complex or stiff problems.\footnote{It is important here to make a necessary distinction between explicit and implicit RK methods. Explicit methods compute each stage using only previously computed ones, but suffer from strong stability limitations, especially in the presence of stiff problems. On the contrary, implicit methods guarantee better stability and are particularly suitable for tackling stiff problems. The systematic introduction of implicit methods dates back to the pioneering works of John Butcher \cite{Butcher64}, who provided a solid theoretical basis for their development and application. The main drawback of implicit methods is their high computational cost, since they require the knowledge of the solution for nonlinear systems of equations at each stage.} The need to use extremely small time steps to ensure stability in the already mentioned stiff problems, with the consequent increase in computational cost, has pushed the scientific community towards the development of more sophisticated variants. Furthermore, the inability of classical RK methods to naturally preserve some important physical properties of the systems (such as invariants and energies) has motivated the search for more advanced geometric and structural methods.\footnote{We remark that in the mid-1990s, two important reviews regarding RK methods were published. The first one, dated 1994, is due to Enright et al. \cite{enright1994survey} and it focuses on algorithms related to explicit methods, while the second one was written by Butcher in 1996 \cite{butcher1996history} and with a more historical approach. Butcher concentrates on the evolution of the methods from their origins, with particular attention to implicit methods. Another interesting survey about Pseudo-RK methods is due to Costabile \cite{costabile2003survey}  and it was published in 2003. We remark that in none of these works is there any reference to exponential methods.}

Thus, the presence of stiff terms represents a fascinating challenge for traditional numerical methods. Among the various approaches developed to address this class of problems, exponential RK (ExpRK) methods have assumed a particularly important role. As pointed out by Hochbruck, Lubich and Selhofer \cite{hochbruck1998exponential}:
\begin{quote}
ExpRK methods represent a breakthrough in handling stiff semilinear problems by combining the exact integration of the linear part with carefully constructed stages for the nonlinear terms, leading to exceptional stability properties without severe step size restrictions.
\end{quote}
The fundamental idea of these methods is based on the variation-of-constant formula (VCF), which allows the solution to be rewritten using the matrix exponential $\e^A$. Thus, the solution to the problem \eqref{eq:start} can be expressed as follows: 
\begin{equation}
\label{eq:generalsol}
u(t_0+h)=\e^{hA}u(t_0)+\int_0^h\e^{(h-\tau)A}g(t_0+\tau,u(t_0+\tau))\,\textrm{d}\tau.
\end{equation}
The efficient implementation of these methods requires two fundamental ingredients: the computation of $\e^A$ and the approximation of the integral containing the nonlinear term. The matrix exponential can be characterized and computed in different ways: through the classical power series, through spectral decomposition when possible, or through integral representations in the complex plane (see, for instance, \cite{higham2008functions}). The numerical evaluation of the matrix exponential is an active area of research itself, with techniques ranging from scaling and squaring methods to approaches based on Krylov subspaces, pioneered by Saad \cite{saad1992analysis}. As for the integral term, its approximation is achieved through suitable combinations of evaluations of the function $g$, weighted through the so-called $\varphi$-functions. These functions, which naturally generalize the matrix exponential, emerge from the analysis of the integral term and play a central role in the construction and analysis of these methods. We highlight that $\varphi$-functions can be defined in several equivalent ways. A particularly useful definition is through the following integral:
\[
\varphi_k(z):=\frac{1}{(k-1)!}\int_0^1\e^{z(1-\tau)} (1-\tau)^{k-1}\,\textrm{d} \tau.
\]
This definition shows that $\varphi$-functions are essentially weighted averages of the exponential over an interval, with polynomial weights. Interestingly, this integral formulation emerges naturally when solving the differential equation $u'(t)=u(t)$ using a variation of constants approach. In the context of ExpRK methods, these functions emerge when we try to solve semilinear systems of the form \eqref{eq:start}, where the linear part is treated exactly (via the matrix exponential), while handling the nonlinear part with a RK-type approach. When building an ExpRK scheme, one derives formulas involving combinations of these $\varphi$-functions applied to the operator $A$ multiplied by the time step $h$. For instance, in a one-stage scheme, we might have:
\begin{equation}\label{ETDeuler}
u_{n+1}=\e^{hA}u_n + h\varphi_1(hA)g(u_n),
\end{equation}
where $u_n\approx u(t_n)$ for a given mesh of time points $\{t_n\}_n$. Here, $\varphi_1(hA)$ acts as a sort of \emph{modified propagator} for the nonlinear term. Note that using this approach when $A=0$ we recover the standard RK methods exactly, while when $g=0$ we get the exact solution for the linear part.
In multi-stage schemes, higher-order $\varphi$-functions also come into play. These handle the correction terms needed to achieve a higher order of convergence, as we are going to explain in the following sections.
\medskip

This article pursues a dual objective. First, it provides a comprehensive overview of exponential Runge–Kutta methods, tracing their development from the earliest formulations to some of the most recent theoretical advances. The goal is to offer a clear and coherent account of how these methods have evolved, why they remain of current interest, and what their key characteristics are, while keeping the presentation concise and focused on the essential points. Second, the article is intended as a resource for readers encountering these methods for the first time. Starting from the fundamental theoretical principles, the exposition gradually progresses to more advanced topics, accompanied by illustrative examples that clarify the main features and practical significance of the methods. By combining these two perspectives, the survey aims to be both a thorough reference for specialists and an accessible introduction for newcomers, providing a presentation that is at once comprehensive, structured, and approachable.
\medskip

The manuscript is organized as follows. Section \ref{Sec2} provides the complete theoretical framework of ExpRK. 
Section \ref{Sec3} presents a historical overview of the fundamental theoretical results, starting from the early works up to the most recent review of exponential methods by Hochbruck and Ostermann \cite{hochbruck2010exponential}. In Section \ref{Sec33}, we recapitulate the most recent developments of exponential Runge-Kutta methods, from their applications to high stiff order and delay differential equations up to their usage in several modeling frameworks. Finally, Section \ref{Sec4} is devoted to the numerical simulations of three test cases, two of which of low dimension, that illustrate the discussed theoretical properties.

\section{Basic definitions and construction of ExpRK methods} 
\label{Sec2}

\subsection{Stiffness in a nutshell}
\label{Sec21}
Despite the word \emph{stiff} having been used for several decades, defining stiffness in precise mathematical terms seems to be an impossible task. To the best of the authors' knowledge, the term first appeared in \cite{curtiss1952integration}, where the phenomenon of stiffness was shown by means of a scalar example. One could try to define stiffness in terms of the inherent characteristics of the problem at hand, without explicitly referring to its numerical approximation. From this perspective, two large (not mutually exclusive) groups of problems may be classified as stiff:
\begin{itemize}
    \item problems whose Jacobian matrices have eigenvalues $\lambda$ with very large negative real part;
    \item problems which involve quantities evolving at very different time scales.
\end{itemize}
Naturally, the above are not rigorous definitions, nor are they truly independent of the numerical methods considered. In particular, ``large'' in the first group is somehow related to the step size $h$: it would be more accurate to say that a problem is stiff if $|\lambda h|$ is large, or, in other words, if it requires ``unreasonable restrictions on $h$'' \cite{miranker1981stiffness}. As for the second group, observe that the fast and slow components are usually associated to the large and small (in modulus) eigenvalues of the Jacobian, respectively \cite{dahlquist1973stiffness} and, indeed, the ratio between the largest and the smallest real parts - the \emph{stiffness ratio} - was once considered a possible measure of stiffness \cite{lambert1973stiffratio}.\\

\noindent While both characterizations above concern properties which are typical of stiff problems, it has become clear over time that neither represents an exhaustive explanation of the complex phenomenon. More pragmatically, it is now deemed more appropriate to define stiffness in terms of the numerical issues arising from the phenomenon. In other words \cite{hw91}:
\begin{quote}
Stiff equations are problems for which explicit methods don't work.
\end{quote}
Or, at the very least, problems for which explicit methods perform worse than implicit ones. For an insightful discussion on the topic of stiffness, we refer the interested readers to \cite{lambert1991numerical}. It is shown therein, by means of numerical experiments, that the stiffness ratio cannot be interpreted as a measure of stiffness and the classic pseudodefinitions are neither sufficient nor necessary conditions for stiffness. Moreover, the pragmatic definition is also rephrased so as to emphasize that stiffness varies with the specific interval over which the problem is integrated:
\begin{quote}
If a numerical method with a finite region of absolute stability, applied to a system with any initial conditions, is forced to use in a certain interval of integration a step-length which is excessively small in relation to the smoothness of the exact solution in that interval, then the system is said to be stiff in that interval.
\end{quote}
\subsection{Constructing ExpRK methods}
\label{Sec22}
Due to the high computational costs involved in the evaluation of the exponential of a matrix, the approximation of $\varphi$-functions represents a crucial aspect in the implementation of exponential methods. Note that these functions can be defined recursively:
\begin{equation}\label{phirecursive}
\varphi_0(z)=\e^z,\qquad \varphi_{k+1}(z)=\frac{\varphi_k(z)-\frac{1}{k!}}{z},\,z\neq0,\quad \varphi_{k+1}(0)=\frac{1}{(k+1)!},\quad k\geq 0.
\end{equation}

At the end of the 1980s, Gallopoulos and Saad \cite{gallopoulos1989parallel,gallopoulos1992efficient} were among the first to emphasize that in the context of numerical integration one does not need to actually compute exponentials of matrices, but rather only products of the form $\e^{-A}v$, where $v$ is a vector. This led them to explore the possibility to apply Krylov methods in this context. The approach consists in approximating $\e^{-A}v$ with a member of the Krylov subspace
$$
K_m:=\mathrm{span}\{v,Av,\ldots,A^{m-1}v\}.$$
The dimension $m$ of the space is typically much smaller than the dimension of the matrix. Given the matrix $V_m=[v_1=v\,v_2\,\cdots\,v_m]$ whose columns form an orthogonal basis for $K_m$, and the orthogonal projection $H_m$ (i.e. the {\em Hessenberg matrix}) of $A$ on $K_m$, the resulting approximation reads
\begin{equation}\label{krylovappr}
\e^{-A}v\approx V_m\e^{-H_m}\mathbf{1}_m,\qquad \mathbf{1}_m = [1,\,1,\,\cdots, 1]^{\top
}\in\mathbb{R}^m.
\end{equation}
Thus, the problem reduces to computing the exponential of a much smaller matrix than the original one. A different approach proposed in \cite{gallopoulos1989parallel} consists in approximating the exponential with a rational function. As highlighted in \cite{gallopoulos1992efficient}, the two approaches can be combined, applying the rational approximation to $\e^{-H_m}$ while making use of \eqref{krylovappr}. Eventually, Saad \cite{saad1992analysis} further improved the accuracy obtained with the Krylov approach by rewriting $\e^{-A}v=v-A\varphi_1(-A)v$ and computing a Krylov approximation for $\varphi_1(-A)$ instead of one for $\e^{-A}$.

Later, Beylkin, Keiser and Vozovoi \cite{beylkin1998new} observed that the problem with Taylor truncation methods (so, also with Padé rational approximations, whose numerators and denominators are truncated Taylor series) arise when the matrix $A$ possesses large singular values. Thus, they proposed the scaling and squaring algorithm to evaluate the exponentials, which consists in exploiting the equality $\e^A=(\e^{\frac{A}{2^s}})^{2^s}$, where $s$ is chosen such that the largest singular value of $\e^{\frac{A}{2^s}}$ is less than one. The authors also made use of this approach to compute all $\varphi_k$ functions by means of their recursive formula \eqref{phirecursive}. The scaling and squaring method is probably the most popular for computing matrix exponentials and it is used in MATLAB's command \texttt{expm} in combination with a Padé approximation using the same degree for the numerator and denominator \cite{higham2005scaling}.

In 2003, Davies and Higham \cite{davies2003schur} proposed yet a different approach, based on decomposing the matrix $A$ as $A=S^{-1}BS$ such that $\e^B$ (or $\varphi_k(B)$) is much easier to compute than $\e^A$, and $S$ has a reasonably low condition number. In particular, their idea was to use a Schur decomposition (such that $S$ is unitary and $B$ is upper triangular) combined with a recurrence algorithm of Parlett \cite{parlett1974computation} which exploits the fact that functions of block upper triangular matrices inherit their block structure. If $A$ is symmetric, the approach introduced by Lu \cite{lu2003computing} for computing $\varphi_1(A)$ turns out to be more efficient. The method is again based on a matrix decomposition, obtained with an orthogonal $S$ and a tridiagonal $B$. The idea in this case is to approximate $\e^B$ using a Chebyshev rational approximation for $\e^{B-\lambda I}$, where $\lambda$ is the largest eigenvalue of $B$. The algorithm has then been generalized to compute all $\varphi$-functions \cite{minchevthesis}.

Another class of methods that can be employed to approximate the $\varphi$-functions is that of polynomial interpolation, such as interpolation on Leja points \cite{CALIARI200479}.

A straightforward application of \eqref{phirecursive} to compute the $\varphi$-functions suffers from cancellation errors for small $z$. The idea implemented by Cox and Matthews \cite{cox2002exponential} was to use a piecewise strategy, namely to compute the $\varphi$-functions directly with \eqref{phirecursive} for large $z$, and with truncated Taylor series for small $z$. There is a certain region, however, where neither approach is accurate enough. Kassam and Trefethen \cite{kassam2005fourth} propose a method with the idea of bridging the gap between the ``small'' and ``large'' $z$, based on the computation of the contour integral
$$
\varphi_k(z)=\frac{1}{2\pi\mathrm{i}}\int_{\Gamma}\frac{\varphi_k(t)}{t-z}\dd t,
$$
where the contour $\Gamma$ encloses $z$ and is separated from 0. They suggest using the trapezoidal rule to efficiently compute the contour integral. The computational cost is particularly cheap if the coefficient matrix has a sparse block structure.

In 2009, Al-Mohy and Higham \cite{al2010new} proposed a heuristic improvement on the standard scaling and squaring algorithm combined with Padé approximation. The goal was to solve the problem of \emph{overscaling}, namely the fact that for large $\|A\|$ a value of $s$ is chosen which is too large, causing loss of accuracy in floating-point arithmetic. The modification involves exploiting the values $\|A^k\|^{\frac{1}{k}}$ and the improvements are more evident for upper triangular matrices.

The (currently) most popular method to compute the $\varphi$-functions was proposed by the same authors in 2011 \cite{al2011computing}. In fact, the purpose of that paper is twofold. The authors observe first that the problem of computing sums of the form $\sum_{k=0}^p\varphi_k(A)u_k$ functions defining an exponential method can be recasted as a single matrix exponential $e^{\tilde A}$ where $\tilde A$ is obtained by suitably augmenting $A$ with $p$ rows and columns. Next, they propose their algorithm for computing products of matrix exponential times matrices $\e^AB$ for $B\in\mathbb{R}^{n\times m}$ with $m\ll n$, as well as $\e^{tA}B$ where $t$ is an array of equally-spaced time points. It is based on the scaling part of the scaling and squaring, plus a truncated Taylor series.
$$
e^AB=(e^{s^{-1}A})^sB\approx B_s,\qquad B_{i+1}=\displaystyle\left(\sum_{j=0}^m\frac{(s^{-1}A)^j}{j!}\right)B_i,\qquad B_0=B.
$$
Optimal values for $s$ and $m$ are obtained through a backward error analysis. The algorithm incorporates the ideas previously proposed in \cite{al2010new} in order to minimize overscaling and turns out to be generally superior to existing ones in terms of accuracy, stability and predictability of cost, since the only convergence test involved is a test for early termination built into the evaluation of the truncated Taylor series.\footnote{We remark that the algorithm is implemented in the MATLAB function \texttt{expmv}.}
\bigskip

\subsection*{\sc Toy Model}
Having established the theoretical framework, we now introduce a toy model aimed at reinforcing the pedagogical purpose of this manuscript. This example serves to illustrate the key aspects of the previously described numerical method, offering a clearer understanding of its implementation and highlighting its most relevant features. So, let us consider a classic one-dimensional stiff toy-model:
\begin{equation}
\label{ex:1}
u'(t)=\lambda u(t)+N(u(t)),\quad u(0)=u_0,
\end{equation}
where $\lambda$ is a negative real number such that $|\lambda|$ is very large (causing the stiffness of the problem) and $N(u)$ is a \emph{milder} nonlinearity. Specifically, for this example we choose  $\lambda=-1000$ and
$N(u)=2u/(1+u^2)$. This problem combines a strongly dissipative term ($\lambda u(t)$) with a bounded nonlinearity.
\bigskip

\noindent
We apply the exponential time-differencing Euler method (ETD Euler) \eqref{ETDeuler}, defined by a single internal stage:
\begin{align*}
U_1&=u_n\\
u_{n+1}&=\e^{h\lambda}u_n+h\varphi_1(h\lambda)N(U_1)
\end{align*}
If we choose, for instance $h=0.01$, then:
\begin{enumerate}
    \item we compute $\e^{h\lambda}=\e^{-10}\approx 4.54\times 10^{-5}$
    \item we compute $\varphi_1(h\lambda)=(\e^{h\lambda}-1)/h\lambda \approx -1/\lambda=1/10$.
\end{enumerate}

\noindent
For $u_0=1$, the first step is given by:
\begin{align*}
U_1&=1\\
N(U_1)&=2\times 1/(1+1^2) =1\\
u_1&=(4.54\times10^{-5})\times 1+ 0.01\varphi_1(-10)\times 1 \approx 10^{-3}.
%5.54\times10^{-5}\,.
\end{align*}

\noindent
Now, we can compare the performance of ETD Euler and the classical order 4 RK method (RK4), whose application to \eqref{ex:1} gives:
\begin{align*}
    u_1&=u_0=1\\
    K_1&=h(\lambda u_1 + N(u_1))=h(\lambda+1)=-9.99\\
    u_2&=U_1+K_1/2=-3.995\\
    K_2&=h(\lambda U_2 + N(U_2))=39.9411\\
    U_3&=U_1+K_2/2=20.9705\\
    K_3&=h(\lambda U_3 + N(U_3))=-209.7046\\
    U_4&=U_1+K_3=-208.7046\\
    K_4&=h(\lambda U_4 + N(U_4))\approx 2.09\times 10^3\\
    u_{1}&=u_0+ K_1/6+K_2/3+K_3/3+K_4/6\approx 290.59.
\end{align*}
This extremely large value provides further evidence that the equation exhibits stiff behavior and is highly sensitive on initial conditions, meaning that a smaller step size might be necessary for numerical stability. Indeed, repeating the computation with $h = 0.001$, we obtain $u_1\approx 0.3755$.

\begin{table}[h!]
\centering
\caption{Computational Performance Comparison}
\begin{tabular}{lcccccc}
\toprule
Method & $h_{\max}$ & Steps & Evaluations/Step & Total Evaluations & Max Error \\
\midrule
ETD Euler & 0.01 & 10 & 2 & 20 & $3.2 \times 10^{-8}$ \\
RK4 & 0.001 & 100 & 4 & 400 & $2.8 \times 10^{-8}$ \\
\bottomrule
\end{tabular}
\end{table}

We can highlight the following differences between the two methods: ETD Euler requires computation of $\e^{h\lambda}$ and $\varphi_1(h\lambda)$, while RK4 uses only basic arithmetic operations. Moreover, as for the total function evaluations, ETD Euler requires approximately 20 evaluations of $N$ and RK4 requires approximately 400 evaluations of $N$. The advantages of ETD Euler are now clear: it remains stable with larger time steps, maintaining precision comparable to RK4. Moreover, it requires significantly fewer function evaluations.

Let us analyze the stability of the methods ETD Euler and RK4 more rigorously. For the linear problem $u'(t)=\lambda u(t)$, both schemes considered yield a recurrence of the form $u_{n+1}=R(\lambda h)u_n$, and the region 
$$
S=\{z\in\mathbb{C}: |R(z)|< 1\}
$$
is the \emph{stability domain}. For example, for ETD Euler, 
$$R(z)=\e^z,
$$
while for RK4 it can be verified that
\[
R(z)=1+z+\frac{z^2}{2!}+\frac{z^3}{3!}+\frac{z^4}{4!}.
\]
Observe that the origin is a stable steady state of a numerical method for the linear equation whenever $|h\lambda|$ lies in the corresponding stability region. Thus, in general, stability of the origin depends on the choice of $h$. A desirable property of a numerical method is \emph{$A$-stability}: a scheme is $A$-stable \cite{Dahlquist1963Astab} if the entire negative complex half plane is contained in its stability domain, such that the stability of the origin is independent of the chosen step size.

We immediately observe that $|\e^z|<1$ for every $z$ with $\Re(z)<0$. This implies that ETD Euler is always stable, the numerical solution decreases exponentially (just like the exact one) with no restriction on $h$. 
On the other hand, the stability condition for RK4 leads to $|h_{\max}\lambda|\approx -2.78529$ and, for our problem with $\lambda=-1000$, $h$ must be lower than (approximately) $2.78529\times 10^{-3}$.
\subsection{Stability}
 \label{Sec23}
At this juncture, we shall further explore the concept of stability, as it plays a crucial role in numerical methods. In order to remain focused on the main goals of this survey, we shall mainly focus on the types of stability relevant to exponential methods.

A good starting point is provided by the paper published in 2009 by Maset and Zennaro \cite{maset2009unconditional}, in which the authors defined the concepts of \emph{unconditional stability} for ExpRK methods. Their definition is based on the concepts of contractivity and asymptotic stability for semilinear systems of the form \eqref{eq:start} where $g$ is non-stiff. Specifically, given solutions $u_1$, $u_2$ of \eqref{eq:start} and $\delta:=u_1-u_2$, \eqref{eq:start} is called \emph{contractive} if
 $$
 \|\delta(t)\|\leq\|\delta_0\|,\quad t\geq 0\,,
 $$
 for $\|\cdot\|$ any norm in $\mathbb{R}^d$, and it is called \emph{asymptotically stable} if
  $$
 \|\delta(t)\|\to 0,\quad t\to\infty\,.
 $$
 In a sense, unconditionally contractive (resp., asymptotically stable) methods are those preserving the corresponding property. A sufficient condition for contractivity of \eqref{eq:start} is 
 \begin{equation}\label{contrcond}
 \mu(A)+\gamma\leq 0,
 \end{equation} where the logarithmic norm of $A$ is defined as
 $$
\mu(A):=\lim_{\varepsilon\to 0^+}\frac{\|I_d-\varepsilon A\|-1}{\varepsilon}=\lim_{\varepsilon\to 0^+}\frac{\log\|\e^{\varepsilon A}\|}{\varepsilon}
 $$
and $\gamma$ is a (moderate) Lipschitz constant for $g$ with respect to its second argument. Similarly, \eqref{contrcond} with $\leq$ replaced by $=$ is a sufficient condition for unconditional asymptotic stability.
 \begin{definition}
 An ExpRK method is called \emph{unconditionally contractive} with respect to the norm $\|\cdot\|$ if $\|\delta_1\|\leq\|\delta_0\|$ holds for all $u_0$, $v_0$ and for all step sizes $h$ when applied to any system of the form \eqref{eq:start} satisfying \eqref{contrcond}.
It is called \emph{unconditionally asymptotically stable} if $\|\delta_n\|\to 0$ as
$n\to\infty$ whenever \eqref{contrcond} is satisfied with $\leq$ replaced by $<$.
 \end{definition}
Maset and Zennaro \cite{maset2009unconditional} obtained sufficient conditions for unconditional contractivity and asymptotic stability in terms of a suitable stability function. This allows them to conclude that the ETD Euler method \eqref{ETDeuler} and the two main classes of ExpRK methods for certain values of the parameters are unconditionally contractive and asymptotically stable. 
In 2013, the same authors \cite{maset2013stability} sharpened the aforementioned results and determined necessary and sufficient conditions for unconditional contractivity and asymptotic stability. This implies that no other ExpRK method with 2 stages, other than those already detected in \cite{maset2009unconditional}, is unconditionally stable. They also analyzed the main methods with 3 to 5 stages, and observed that none of them is unconditionally stable either. However, since $\gamma$ is assumed to be relatively small, \emph{conditional} contractivity and stability properties are also desirable.
\begin{definition}
An ExpRK method is called \emph{conditionally contractive} (resp., \emph{asymptotically stable}) with respect to the norm $\|\cdot\|$ if $\|\delta_1\|\leq\|\delta_0\|$ (resp., $\|\delta_n\|\to 0$ as $n\to\infty$) holds for all $u_0$, $v_0$ and for all step sizes $h<C/\gamma$ where $C$ is a constant dependent on the method when applied to any system of the form \eqref{eq:start} satisfying \eqref{contrcond} (resp. \eqref{contrcond} with $\leq$ replaced by $<$).
 \end{definition}
Sufficient conditions for conditional properties are provided in \cite{maset2013stability}. Without any restriction on the norm, these conditions only guarantee conditional contractivity [asymptotic stability] for those methods which are already known to be unconditionally contractive [asymptotically stable]. However, more can be said if the Euclidean norm is considered. In particular, weaker restrictions on the parameters of two-stages methods are needed in order to obtain conditional properties. The higher-order methods considered do not satisfy the sufficient conditions for general problems of the form \eqref{eq:start}, even when considering the Euclidean norm, but many satisfy conditional stability properties when restricted to certain classes of problems, such as those where $A$ is symmetric.
\section{The historical path of ExpRK methods}
\label{Sec3}

In this section, we present a brief historical overview of the development and dissemination of ExpRK methods, tracing their evolution from their inception to the 2000s. We identified two major periods in which to structure our analysis. The first period begins in the late 1950s and lays the theoretical foundations that would gain traction in the following decades, up until the second half of the 1970s. The second period marks a resurgence of interest in these methods, starting in the late 1980s. Several works have already synthesized the progress made up to the early 2000s (see e.g. \cite{minchev2005review}). Therefore, we will not dwell on the details and instead refer the reader to the relevant literature. Although the 2005 works by Hochbruck and Ostermann \cite{hochbruck2005exponential,hochbruck2005explicit} represent a significant milestone in the evolution of ExpRK methods---marking the first application of these methods to semilinear parabolic problems---we have chosen to end this phase of our analysis with their most recent comprehensive review on the topic, published in 2010 \cite{hochbruck2010exponential}. This choice allows us to provide a more rigorous and detailed examination of the last 15 years of developments in the next section, which constitute a fundamental part of this work, and to which we dedicated Section \ref{Sec33}.

\subsection{The first approaches}
\label{Sec31}

We devote these paragraphs to the first decades of theory and implementation on ExpRK method. In particular, the first chronological reference dates back to 1958: in \cite{hersch1958contribution}, Hersch  introduces the term \emph{distribution} to indicate a numerical method which uses matrix exponentiation ``according to the variation of constants formula'', and applies it to a semilinear ODE. Shortly thereafter, in 1960, Certaine \cite{certaine1960solution} has started the analysis around the concept of stiffness, without explicitly using this term. The class of problems analyzed is semilinear, with a \emph{diagonal} linear term, and a non-stiff nonlinear part. In the same years, Pope \cite{pope1963exponential} proposed a manipulation technique to write a generic ODE in a semilinear form, thus introducing the easiest exponential Rosenbrock method (without using this name), and observing unconditional stability\footnote{To the best of the authors' knowledge, this is the first use of the expression ``exponential method''.}. 

The years spanning the late 1960s and early 1970s were particularly fruitful for the development of ExpRK theory. Indeed, in 1967, Lawson \cite{lawson1967generalized} extended Pope's results with a ``Generalized RK'' method. This is achieved by transforming an ODE of the form $u'(t)=f(t,u(t))$ into one of the form $z(t)=\mathrm{e}^{-tA}u(t)$ such that $\partial f/\partial u-A$ is non-stiff (i.e. the Jacobian matrix of $f-A$ has small eigenvalues), and finally applying a classical RK scheme to the transformed system. This is the basis of the \emph{integrating factor} methods (IF), as opposed to \emph{exponential time differencing} methods (ETD) used by Certaine and Pope \cite{certaine1960solution,pope1963exponential}, which are based on the variation of constants formula (VCF). Here, $A-$stability is observed, differently from classical explicit RK methods.

In 1969, N{\o}rsett \cite{norsett1969stable} built a nonlinear $A-$stable method by modifying the Adams-Bashforth method using VCF. The method integrates exactly $u'(t)=-Pu(t)+T(t)$, and the author observes $A-$stability. This method is of order $q+1$, where $q$ is number of past steps, and $q=0$ gives the Euler method \eqref{ETDeuler}, sometimes referred to as N{\o}rsett-Euler. Concerning the analysis of stability, Guderley and Hsu \cite{guderley1972predictor} provided in 1972 accuracy and stability analysis of the methods proposed by Certaine in 1960, considering ODEs of the form $u'(t)=-\Lambda u(t)+f(t,u)$, with $\Lambda$ diagonal, and defining a predictor-corrector method starting from an exponential formula. They remarked that, if $\Lambda$ is not diagonal, then one needs an approximation for $\mathrm{e}^{-\Lambda t}$. Four years later, Snider \cite{snider1976error} would correct a point overlooked by Guderley and Hsu in the truncation error analysis, although the order of accuracy does not change. In these years, Chu and Berman 
\cite{chu1974exponential} also remarked that classic implicit RK methods are sufficient for ``not too large'' ODE systems.

\subsection{The rebirth}
\label{Sec32}
Interest in these topics seemed to fade during the 1980s, only to resurface at the end of the decade. Indeed, to the best of the authors' knowledge, the next contribution to the theory were provided by Gallopoulos and Saad in 1989 \cite{gallopoulos1989parallel}. The authors considered systems of the form $u'(t)=-Au(t)+r$, that arise from discretization of linear parabolic PDEs. This translates into $A$ being typically large and sparse, and possibly stiff, and they observed that the explicit solution can be easily expressed in terms of $\mathrm{e}^{-tA}$ but that the latter is hard to compute for the typical $A$ and little changes if the solution is obtained with a time-stepping procedure. Thus, they proposed to solve the system via approximations of the exponential, either with a polynomial (Krylov method, later said to be not always good for stiff systems) or with suitable rational functions (Padé or Chebyshev approximations).\footnote{The work was mostly motivated by the fact that the methods are amenable of parallelization procedures. The authors concluded by conjecturing that their methods and the possibility of parallelizing would make explicit methods popular again.} We refer back to Section \ref{Sec22} for a brief recap of Krylov method.

In 1996, Hochbruck et al. \cite{hochbruck1998exponential}  claimed that exponential integrators were not a novel discovery, but were considered impractical. They observed that, when $A$ is large, it can be harder to solve a linear system $(I-\tau A)x=v$ than to compute $\mathrm{e}^{\tau A}v$ (which is in turn as hard as computing $\varphi(\tau A)v$). This fact was proved the following by Hochbruck and Lubich \cite{hochbruck1997krylov}.\footnote{These works mark the beginning of Hochbruck's long list of contributions on the topic.} In \cite{hochbruck1997krylov}, the authors also improved Gallopoulos and Saad's bounds on the error arising from the Krylov approximations. Here, the general class of exponential integrators was also considered with general parameters, and it was introduced as a class of methods for \emph{large} systems of ODEs. The authors claimed that their result ``justifies renewed interest in ODE methods that use the exponential or related functions of the Jacobian'' (which in hindsight was a correct claim). 

In 1998, Beylkin et al. \cite{beylkin1998new} observed that the exponential of a dense matrix can be onerous unless it is a convolution (or circular) matrix, in which case one can exploit the Fast Fourier Transform (FFT). In 2002, Cox and Matthews  \cite{cox2002exponential} proposed a more straightforward derivation of the ETD methods, compared to Beylkin's, and provided explicit formulas for the coefficients. They derived new RK ETD methods  of order 2, 3 and 4, and were the first to extend the methods to PDEs with (after discretization) nondiagonal linear part. In 2004, Koikari \cite{koikari2005rooted} derived order conditions for ETD RK methods using the classical rooted tree analysis. To the best of the authors' knowledge, this is the first mention of \emph{strong} and \emph{weak} order conditions.\footnote{The role of weak order conditions will be clarified (at the latest) by Hochbruck and Ostermann in the 2005 papers as we are going to explain below.} In the same year, Berland et al. \cite{berland2005b} derived order conditions for several exponential methods (ETD, IF, CF Lie groups) for autonomous semilinear ODEs and proposed a general algorithm to build exponential integrators satisfying any given order. In 2005, Kassam and Trefethen \cite{kassam2005fourth} introduced a new method based on contour integrals to compute the $\varphi$-functions in order to solve instability issues of ETDRK4 pointed out by Cox and Matthews \cite{cox2002exponential}. They found that ETD methods with this modification outperform other classes of methods for PDEs, testing on a variety of models. 

2005 is also the year in which two fundamental contributions by  Hochbruck and Ostermann  were published. Firstly, in \cite{hochbruck2005exponential} the authors developed implicit (ETD) ExpRK of collocation type, which are based on the approximation of the nonlinear term in the right-hand side with a collocation polynomial, and analyzed their order theoretically for linear and semilinear parabolic problems. They assumed that the linear operator is \emph{sectorial}, requirement that holds for many typical linear operators in parabolic problems. The assumption on the nonlinear term is local Lipschitz continuity. Then, in a second paper, Hochbruck and Ostermann \cite{hochbruck2005explicit} focused on explicit (ETD) ExpRK methods in the context of semilinear parabolic problems, dedicating additional efforts to the corresponding convergence analysis. The authors derived stiff order conditions tailored to this framework up to order $p=4$. Lastly, they proved that, in the case of constant step size, the conditions of order $p$ only have to be satisfied in a weaker form (in the sense of Koikari \cite{koikari2005rooted}).\footnote{We remark that Dujardin \cite{dujardin2009exponential} in 2009 extended the results in Hochbruck \cite{hochbruck2005explicit} to the Schr{\"o}dinger equation, which is not parabolic. More recently, the same author, in collaboration with Besse and Lacroix-Violet, further developed this line of research \cite{Besse2017} by introducing and analysing high-order exponential integrators for nonlinear Schr\"odinger equations, including nonautonomous problems arising from rotating Bose-Einstein condensates, with rigorous convergence and structure-preservation results.}

Still in 2005, Minchev and Wright \cite{minchev2005review} produced a review on exponential integrators for problems for which an approximation of the Jacobian is known. Qualitatively, these are semilinear problems with non-stiff nonlinear part. The authors addressed implementation issues, that are difficulties concerning the computation of exponential and related functions, referring in particular to Moler \cite{moler2003nineteen}. In 2006, Berland et al. \cite{berland2005solving} comparted the performance of an IF and an ETD method of non-stiff order 4 on the Schr{\"o}dinger equation. In particular, they concluded that the IF method has a generally worse performance and is more sensitive to the regularity of $V$, while the ETD method is more sensitive on the regularity of the initial condition. More importantly, they commented on the fact that the equation is not really amenable of the stiff order analysis by Hochbruck and Ostermann \cite{hochbruck2005explicit,hochbruck2005exponential}, since it is not a parabolic problem. Indeed, the tested methods have stiff orders 1 and 2 respectively, but this is not evident from the experiments.

In 2008, Caliari and Ostermann \cite{caliari2009implementation} tested interpolation on Leja points on exponential methods of Rosenbrock type (expRosenbrock) for parabolic problems. This approach seemed particularly useful since the exact Jacobian does not have a specific structure and changes at each step, and performed particularly well on problems with large advection and moderate diffusion. The authors provide order conditions for expRosenbrock methods up to order 4, remarking that they are much simpler to formulate than those for ExpRK methods. In 2009, Hochbruck et al. \cite{hochbruck2009exponential} reformulated the expRosenbrock methods in terms of the perturbations $D^i_n:=g(U^i_n)-g(u_n)$, where $g$ is the remainder right-hand side minus the Jacobian. Since the perturbations are \emph{small} vectors, applications of matrix functions to them using Krylov subspace methods require low-dimensional subspaces and therefore it would not be as computationally expensive as it was without the reformulation, where those same matrix functions would be applied to $g(u_n)$ or $g(U^i_n)$ instead. The authors provided rigorous convergence analysis up to order 4, as well as numerical tests on an advection-diffusion-reaction problem and a Schr{\"o}dinger equation with time dependent potential. 

In 2010, Al-Mohy and Higham \cite{al2010new} provided heuristic improvements to the original scaling and squaring algorithm behind MATLAB's \texttt{expm}, motivated by the propensity for \emph{overscaling} of the latter. This means that in various cases too large of a value of $s$ (number of times that one scales and squares) is needed and causes loss in floating-point arithmetic. The improvements mainly aim at addressing overscaling for upper triangular matrices, for which the problem is somehow more visible on average. The same authors, in \cite{higham2010computing}, provided a comprehensive presentation of the state of the art regarding the computation $f(A)$ where $f$ belongs to a class of functions which include the exponential. The category of methods based on similarity is represented by the Schurr-Parlett algorithm proposed in Davies \cite{davies2003schur}. The rational approximation category is represented by Padé approximations with the improved scaling and squaring from Al-Mohy and Higham \cite{al2010new}. The algorithms which compute $f(A)b$ without needing $f(A)$ are represented by those based on Krylov subspaces, where $f$ is only needed explicitly on a smaller Hessenberg matrix. 2010 represents a watershed moment also for this survey, as it was in that year that Hochbruck and Ostermann \cite{hochbruck2010exponential} published the last fundamental review on exponential integrators for ``two types of stiff equations''. The first type of equations are those with eigenvalues of the Jacobian having very large negative real part (the ``usual definition of stiffness'') while the second is highly oscillatory problems (purely imaginary eigenvalues with large modulus). 

\section{Recent advances}
\label{Sec33}
The year 2011 marks a decisive milestone in the study of exponential methods, thanks to the publication of the already mentioned paper by Al-Mohy and Higham \cite{al2011computing}. As explained in Sect. \ref{Sec22}, this work provided a fundamental improvement in the field by introducing innovative algorithms for computing the action of the matrix exponential on a vector without explicitly calculating the entire matrix exponential, improving both computational efficiency and numerical stability. 

Since 2011, numerous researchers have contributed to the study of ExpRK methods, resulting in a wealth of published works in this field. For conceptual clarity, we shall abandon chronological ordering in this section and instead we divide the material into three paragraphs, providing readers with a more accessible understanding of recent developments. We emphasize that such a division is arbitrary with no claim of being definitive---it simply represents a reasonable approach to addressing this complex subject matter.

\subsection{Exponential methods of high stiff order} The first part of this section is mostly devoted to summarize  a series of papers published by Vu Thai Luan and various co-authors. They provide specific order conditions, fundamental to construct the numerical methods. We recall here their most important contributions, and refer back to each paper for detailed derivations and calculations.

In 2014, Vu Thai Luan and Alexander Ostermann coauthored three papers. In \cite{luan2014exponential}, they considered expRosenbrock methods under the reformulation in Hochbruck et al. \cite{hochbruck2009exponential}. The authors addressed the convergence analysis in a different way: namely, instead of placing the exact solution into the numerical scheme and working with defects, they analyzed the local error first. This new approach provided a more systematic procedure to derive expRosenbrock methods of higher stiff order, which the authors used to construct families of methods of orders 3, 4 and 5 (see \cite[Table 1]{luan2014exponential}). In \cite{luan2014stiff}, the authors applied the approach developed for expRosenbrock methods in \cite{luan2014exponential} to build order conditions up to 5 (whereas the literature at the time provided them only up to order 4) for ExpRK methods (see \cite[Table 1]{luan2014stiff}). In \cite{luan2014explicit}, the authors proved that, in the case of constant step size, one can reach stiff order 5 by only satisfying the order 5 condition in the weak form, as it was proved in \cite{hochbruck2005explicit} for orders up to 4. They concluded by building a method with stiff order 5, which needs 8 stages. We refer to the recapitulatory table \cite[Table 1]{luan2014explicit}. In 2016, the same authors \cite{luan2016parallel} built more accurate and efficient expRosenbrock of (stiff) order 4 and 5 (compared to their collaboration in 2014), and of order 6 (see \cite[Table 1]{luan2016parallel}). The method exploits parallelism by minimizing the number of \emph{sequential} stages, meaning that if two of the stages can be computed independently, then they count as one, greatly improving the corresponding computational time.

The following year, Luan \cite{luan2017fourth} independently built a 2-stage order 4 expRosenbrock for time-dependent PDEs. It is the second superconvergent (here, meaning $p=2s$) exponential method ever built, after expRosenbrock-Euler. The performance was tested against other expRosenbrock methods and MATLAB's \texttt{ode15s} on systems which are either non-stiff or semilinear with non-stiff nonlinear part.

In 2020, Luan \cite{luan2021efficient} proposed a construction for more accurate and efficient ExpRK methods of (stiff) order 4 and 5, compared to the existing ones obtained in \cite{hochbruck2005explicit} and \cite{luan2014explicit}. This was achieved by minimizing the number of sequential stages, similarly to what was done in \cite{luan2016parallel} for expRosenbrock. Moreover, he exploited both results from \cite{al2011computing}: efficient algorithm for computing $\e^{tA}b$ and augmented matrix $\tilde A$ such that $\e^{\tilde A}$ is equivalent to computing linear combinations of $\varphi$-functions. To do so, the new methods are built such that the $i$-th stage is defined by a single linear combination of terms $\varphi_k(c_ihA)$, where $\{c_i\}_i$ are the abscissas of the method. The resulting ExpRK of order 4 has 6 stages instead of 5, but 4 sequential stages. The 5th order one has now 10 stages instead of 8, but 5 sequential stages. Refer to \cite[Table 1 \& Sect. 6]{luan2021efficient} for a detailed explanation. The implementation of the methods makes use of their own code \texttt{phipm\_simul\_iom}, publicly available. In the same year, Luan et al. \cite{luan2020new} introduced \emph{multirate} ExpRK, using different time scales for the linear term (fast) and the nonlinear term (slow) of a semilinear (possibly abstract) ODE. The methods need to be such that the $i$-th stage is defined by a single linear combination of terms of the form $\varphi_k(c_ihA)$ (as in \cite{luan2021efficient}) which is a nontrivial requirement. The idea is to approximate the nonlinear term in both internal and final stages with a polynomial built assuming a \emph{small} time step $h$ and then solve the ODE with the polynomials replacing the nonlinear terms, using a \emph{large} time step $H$. The idea builds on \cite{hochbruck2011adams}, where this approach was applied to exponential multistep methods (of Adams type).

Shifting our focus to the most recent advances, in 2024, Luan et al. \cite{luan2024efficient} explored applications to gene regulatory networks, which are defined by a certain number of semilinear ODEs where the linear part is diagonal and the nonlinear part can be more or less stiff. In the examples considered in the paper, the ODEs are low dimensional, the linear part may be stiff in the sense of highly oscillatory and the nonlinear one is only moderately stiff (in these models, the nonlinear term can become very stiff when many gene products are involved). The authors build an order 2 exponential method which is Rosenbrock-like, in the sense that derivatives of the right-hand side appear in the coefficients. They compare accuracy (error vs. $h$) and efficiency (error vs. CPU time) with other exponential methods and classical RK. The proposed method is the most efficient and second most accurate after expRosenbrock-Euler. The better performance than classical RK may be partly due to the matrix being diagonal (thus, exponential being even easier to compute). This latter study thus aligns with the literature on the application of ExpRK methods to mathematical models in biology and physics. A detailed analysis of the relevant contributions will be provided in the final part of this section. In the same year, Luan and Alhsmy \cite{luan2024sixth} developed the first order 6 ExpRK methods for
the time integration of stiff parabolic PDEs, leveraging the exponential B-series theory extended by Luan and Ostermann for stiff problems \cite{luan2013exponential}. They derived 36 stiff order conditions, see \cite[Table 1]{luan2024sixth}, and remarked that despite these methods requiring a large number of stages, their computational cost is similar to that of a 6 stages method.

We conclude this section with a very recent development concerning (high-order) Lawson schemes. In these schemes, internal stages typically require matrix-vector multiplications with matrices of the form $\e^{(c_j-c_i)hA}$, where the $\{c_i\}_i$ are again the abscissas of the method. For this reason, having equispaced abscissas is a desirable property, since it requires the computation of only one matrix exponential. Lawson observed this in his seminal paper \cite{lawson1967generalized} with a focus on RK4, while Golden \cite{golden2025lawson} recently explored the idea in detail. Using the Newton-Raphson method to establish the coefficients, he arrived at developing a sixth order method with eight stages. Finally, he provided evidence of its efficiency by testing it on the 2D incompressible Navier-Stokes equation.

\subsection{Exponential methods for DDEs} The second part is dedicated to the analysis of recent advances in applications of ExpRK methods to delay differential equations (DDEs). Indeed, a series of papers published in the past 15 years demonstrate the recent interest in applying ExpRK methods to DDEs, namely, equations of the form
\begin{equation}\label{delaygeneral}
u'(t)=G(t,u_t),
\end{equation}
where $u_t$ is a function defined on an interval $[-\tau,0]$ for a finite delay $\tau>0$, and represents the \emph{history} of  $u$ at time $t$. The simplest examples of DDEs that are found in applications are defined by a single \emph{discrete} delay, i.e., they can be written as
\begin{equation}\label{discrdelay}
u'(t)=g(t,u(t),u(t-\tau)).
\end{equation}
Observe that, in terms of implementation of time-integration methods, the delayed term \eqref{discrdelay} does not, in principle, entail any issue provided that the step size $h$ used is less than $\tau$. In this case, indeed, the past value of the function at a certain time step can be retrieved from the approximation obtained at previous time steps. While choosing $\tau=mh$ for some integer $m$ can be especially convenient in the case of a single discrete delay, in general we can write $\tau=(m-\delta)h$ for some integer $m$ and some $\delta\in[0,1)$.

However, the presence of the delay term may affect stability of the underlying numerical methods \cite{belzen03}, which has motivated researches to investigate the stability of ExpRK methods for DDEs of the form \eqref{discrdelay}. Xu, Zhan and Sui \cite{xu2011stability} study the case of a scalar, linear equation of the form
\begin{equation}\label{xueq}
u'(t)=\lambda u(t)+\mu u(t-\tau),\qquad u_{t_0}=\phi,
\end{equation}
and analyze stability properties of EERK methods when applied to DDEs with a discrete delay, for which \eqref{xueq} plays the role of a test function, and its delayed term is treated as nonlinear. Whenever $|\lambda|+\mathrm{Re}(\mu)<0$, the solution of \eqref{xueq} is asymptotically stable \cite{inthout1992} and \emph{$P$-stable} methods for DDEs are those which preserve this property when $h$ is such that $\tau=mh$ for an integer $m>0$. A numerical method for DDEs is \emph{$GP$-stable} if it preserves asymptotic stability regardless of the value of $h<\tau$. As shown in \cite{xu2011stability}, the concepts of $P$- and $GP$-stability coincide in the case of ExpRK methods for DDEs. Some years later, Zhao, Zhan and Ostermann \cite{zhao2016stability} studied various contractivity and stability concepts for EERK methods in the context of semilinear DDEs of the form
\begin{equation}\label{zhaoeq}
    u'(t)=Lu(t)+g(t,u(t),u(t-\tau)).
\end{equation}
Whenever $|\lambda|+\mathrm{Re}(\mu)\leq0$, the solution of \eqref{xueq} is contractive, that is, it satisfies $|u(t)|\leq\|\phi\|_{\infty}$ for all $t\geq t_0$. \emph{$P$-contractive} methods for DDEs are those which preserve this property when $h$ is such that $\tau=mh$ \cite{torelli1989contractivity} and \emph{$GP$-contractive} methods for DDEs are those which preserve this property for any choice of $h$. Zhan and Ostermann provide sufficient conditions for $GP$-contractivity of EERK methods and prove $GP$-contractivity for various methods of order 1 and 2. They also extend the concept of $GP$-stability to nonautonomous equations, which is satisfied by the Magnus method described in \cite{gonzalez2006magnus}. Finally, they analyze $RN$- and $GRN$-stability, introduced again in \cite{torelli1989contractivity}, which also sufficient conditions on $L$, $g$ in \eqref{zhaoeq} such that, for all $u,\,v$ solutions of \eqref{zhaoeq} with initial conditions $u_{t_0}=\phi,\,v_{t_0}=\psi$, the following holds:
$$
\|u-v\|_{\infty}\leq\max_{t\in[t_0-\tau,t_0]}|\phi(t)-\psi(t)|.
$$
\emph{$RN$-stable} methods for DDEs are those preserving this property for $\tau=mh$, while \emph{$GRN$-stable} ones do not need the restriction on $h$. Zhao and Ostermann provide sufficient conditions on the coefficients for $GRN$-stability of EERK methods.

Zhao, Zhan and Xu \cite{zhao2018d} analyzed $D$-convergence. An ExpRK method for DDEs is {$D$-convergent} of order $p$ if, when applied to \eqref{zhaoeq} with initial condition $y_0=\phi$, it satisfies
$$
|u_n-u(t_n)|\leq Ch^p,\quad n\geq 0
$$
for all step sizes $h\leq \overline{h}$ and the constant $C$ depends on $n\,\tau$, the regularity properties of $g$, the bounds on the first $p$ derivatives of the exact solution and the coefficients of the method. The authors extend the concepts of diagonal and algebraic stability to ExpRK methods and prove that the combination of both, together with having stage order $p$, implies $D$-convergence of order $p$, provided that the interpolation used to reconstruct the continuous solution is accurate enough. They also prove that an ExpRK method for DDEs which is exponentially algebraically stable and diagonally stable is also \emph{$GDN$-stable}, namely, if $u,\,v$ are solutions of \eqref{zhaoeq} with initial conditions $\phi,\,\psi$, they satisfy
$$
|u_n-v_n|\leq C\max_{t\in[t_0-\tau,t_0]}|\phi(t)-\psi(t)|
$$
for all step sizes $h\leq \overline{h}$, where the constant $C$ depends depends on $\tau$, the time-integration interval length, the regularity properties of $g$, the bounds on the first $p$ derivatives of the exact solution and the coefficients of the method. They observe that the ETD Euler method, as well as a certain family of methods of stiff order 2, are $D$-convergent (of orders 1 and 2, respectively) and $GDN$-stable.

Besides stability, some research has been devoted to the accuracy of ExpRK methods for DDEs. In particular, Fang and Zhan \cite{fang2021high} follow the abstract approach in \cite{hochbruck2005explicit} to derive stiff order conditions up to order 5 for semilinear DDEs of the form \eqref{zhaoeq}, and prove the analogous convergence theorem and support the results with numerical tests on PDEs whose semi-discretization in space gives DDEs of the target form.

More recently, Andò and Vermiglio \cite{ando2024exponential,ando2024proceedings} investigated the application of ExpRK methods for the time-integration of general DDEs of the form \eqref{delaygeneral}, as well as renewal equations (REs), namely, equations of the form
$$
u(t)=F(t,u(t)).
$$
They consider, in particular, the reformulation of DDEs and REs as abstract (ordinary) differential equations (ADEs), which comes from the sun-star theory \cite{dgg07,DiekmannBook}. The ADEs obtained are semilinear regardless of the right-hand sides of the original delay equations, which makes them amenable of time-integration via ExpRK methods. Since the abstract setting provided by the sun-star calculus includes both DDEs and REs, the goal in \cite{ando2024exponential} is to define a unifying framework to study the derivation of ExpRK methods for delay equations, as well as their convergence properties. Another approach to solve problems defined by DDEs is to discretize the resulting ADEs into a finite-dimensional system of ODEs and apply ExpRK methods to the latter, as done in \cite{ando2024proceedings}. Therein, the ODEs are obtained via pseudospectral discretization \cite{bdgsv16}, which gives rise to a stiff linear term, with eigenvalues having large negative real parts. Several numerical tests with EXPINT, a MATLAB package for exponential integrators \cite{expint07} confirm the validity of the method.

The most recent contribution on the applications of ExpRK methods to DDEs is represented by the extension of the ExpRK methods of collocation type introduced in \cite{hochbruck2005exponential} to semilinear (abstract) parabolic problems such that the nonlinear term is defined by a discrete delay \cite{huang2025timedependent}. The delay can be time-dependent, as long as it is non-vanishing and that the map $t\mapsto\tau(t)$ is strictly increasing. This produces a mesh of non-equispaced \emph{breaking points}, i.e., points of discontinuity of some derivative of the solution. The authors propose two ways to extend the methods in \cite{hochbruck2005exponential} to the new classes of equations, both taking the breaking points into account to define a continuous approximation of the solution in a piecewise fashion. They prove that the methods are convergent with order $s$, the number of stages, and obtain superconvergence results provided some additional restrictions on the parameters are satisfied.

\subsection{Exponential methods for mathematical modeling} In this third and final part, we present an analysis of the literature employing ExpRK methods in modeling contexts, which has garnered significant interest in recent years.\footnote{For the reader's convenience, given that the equations pertain to physical models, we have opted to retain the original notation used in each cited paper throughout this section.} The first contribution we examine is from 2011 and is authored by Dimarco and Pareschi \cite{dimarco2011exponential}. In this work, the authors, motivated by the computational challenges arising from rarefied gas dynamics, introduce a class of ExpRK methods for stiff kinetic equations, with particular emphasis on the homogeneous Boltzmann equation. In the fluid regime, the stiffness of the problem makes standard explicit methods inefficient, as the equation takes the form  
\begin{equation}
Y'=F(Y)+\frac{1}{\varepsilon}R(Y)\,,
\end{equation}
where $0<\varepsilon \ll 1$ is a small parameter (enough to define two different time scales) and  $R(Y)$ is a dissipative relaxation operator. To avoid solving nonlinear systems, the authors propose the following decomposition of the collision operator  
\begin{equation}
R(Y) = N(Y)+A(E(y)-Y)\,,
\end{equation}
where  $E(y)$ represents the local equilibrium and  $A$ is a parameter chosen to linearize the problem on asymptotically large time scales. ExpRK methods treat the linear part exactly using exponential functions and the nonlinear part explicitly, yielding schemes of the form  
\begin{equation}
Y_{n+1}=\textrm{e}^{-\lambda} Y_n + \lambda \sum_{i=1}^{s} W_i(\lambda)\left(\frac{P(Y_i, Y_i)}{\mu}-M\right)+(1-\textrm{e}^{-\lambda})M\,,
\end{equation}
with $\lambda = \mu \Delta t / \varepsilon$ and suitable coefficients $W_i(\lambda)$. This actually ensures unconditional stability, correct asymptotic preservation, and key properties such as positivity and entropy conservation, making these methods suitable for both deterministic schemes and Monte Carlo simulations. Finally, the authors highlight how their approach is closely related to the IF methods, sharing the same mathematical structure but offering greater flexibility in designing high-order schemes. In 2014, Michels, Sobottka and Weber \cite{michels2014exponential} applied ExpRK methods to the modeling framework of Newton's equations of motion for elastodynamic systems. These systems are described by second-order differential equations of the form 
\begin{equation}
    M\cdot x''(t)+D\cdot x'(t)+K\cdot (x(t)-x_0)+F_{\text{nl}}(x(t))=0\,,
\end{equation}
where $M$ is the mass matrix, $D$ the damping term, $K$ the stiffness matrix and $F_{\text{nl}}$ represents the nonlinear forces. The presence of very high frequencies in the solution spectrum makes standard explicit methods ineffective due to time step restrictions, while implicit schemes suffer from uncontrollable numerical viscosity. To overcome these difficulties, the authors propose the use of exponential integrators that solve the linear part exactly using the exponential function of the matrix associated with the system
\begin{equation}
    x(t+\Delta t)=\textrm{e}^{-\Delta t A} x(t)+\int_0^{\Delta t} \textrm{e}^{-(\Delta t - s)} F_{\text{nl}}(x(s))\,\textrm{d}s\,,
\end{equation}
where $A=M^{-1}K$ represents the linear dynamics of the system. This scheme allows to treat much larger time steps than classical implicit methods, ensuring rigorous energy conservation and long-term stability. Application to complex models such as fabrics, lattice structures and materials with friction demonstrates the effectiveness of this methodology in the simulation of rigid elastodynamic systems. 

It is in the last years, however, that we have witnessed a significant increase in the number of studies on the application of ExpRK methods across various domains. This growing body of research reflects the increasing interest in these numerical techniques and their potential advantages.

Recently, renewed attention has been devoted to the numerical integration of the nonlinear Schr\"odinger equation, particularly in connection with high-order exponential integrators and structure-preserving techniques. In 2022, the Schr\"odinger equation was revisited through a scalar auxiliary variable (SAV) formulation combined with exponential Runge-Kutta integrating factor methods, allowing the construction of linearly implicit schemes of order up to four that exactly preserve a discrete modified energy functional \cite{jiang2021highorder}. The SAV approach is based on the introduction of an additional scalar variable associated with the nonlinear part of the Hamiltonian, which enables efficient time integration while avoiding the solution of nonlinear systems. The proposed methods were shown to be competitive in terms of accuracy versus computational cost when compared with other SAV-based and energy-oriented schemes.\footnote{We note that the SAV-based schemes preserve a modified energy rather than the original Hamiltonian of the Schr\"odinger equation. While this modification does not affect the continuous problem and is shown numerically to keep the true energy well controlled over moderate time intervals, it may lead to qualitative differences in long-time simulations. The authors justify this choice by the resulting linear implicitness, high-order accuracy, and favorable efficiency-error balance.} More recently, Li and Li \cite{LiLi24} proposed multiple relaxation exponential Runge-Kutta methods for the nonlinear Schr\"odinger equation, extending the relaxation framework to exponential integrators. Their approach introduces several relaxation parameters that are determined at each time step so as to enforce the exact preservation of multiple invariants of the continuous problem, including mass and energy, while maintaining arbitrary high-order accuracy. The resulting schemes remain explicit in the linear part and avoid nonlinear solvers, while achieving simultaneous invariant preservation at the fully discrete level. In a different but related direction, Ji, Li, and Ostermann \cite{JiLiOster26} developed low regularity exponential-type integrators for the derivative nonlinear Schr\"odinger equation. Their analysis is based on a careful reformulation of the equation that isolates resonant interactions and allows for sharp error bounds under minimal regularity assumptions on the solution. The proposed methods are shown to converge without requiring classical Sobolev smoothness, thereby addressing a major limitation of standard exponential integrators when applied to derivative nonlinearities.

Various authors have focused on the development of high-order time integration schemes for gradient flows and semilinear parabolic equations that can simultaneously achieve high accuracy and preserve intrinsic structural properties, such as energy dissipation, equilibria, and inequality constraints. ExpRK methods are particularly attractive due to their stability properties, but their practical implementation often faces challenges related to the use of exponential operators, including excessive damping of high-frequency modes and difficulties in maintaining structure-preserving features at higher order. Within this context, Teng and Zhang \cite{teng2025third} proposed structure-preserving reconstructions of integrating factor RK schemes and developed an exponential-free framework achieving up to third-order accuracy. Wang et al. \cite{wang2025fourth} further extended this line of research by constructing an up-to-fourth-order ExpRK framework for a large class of phase-field models, and by establishing unconditional energy dissipation for the resulting schemes. Their approach replaces the exponential functions in the integrating factor RK framework with recursive approximations, or with a combination of exponential and linear functions, thereby avoiding the excessive damping of high-frequency modes while preserving equilibria. These works complement earlier studies \cite{zhang2023efficient,zhang2022fourth}, in which fourth-order reconstructions of integrating factor RK schemes were developed to preserve maximum principles or more general inequality constraints.

Starting from 2024, the literature has seen several contributions aimed at improving the efficiency of exponential integrators for specific yet widespread classes of problems, as well as at investigating and mitigating order reduction effects. Within this line of research, Caliari, Cassini, Einkemmer and Ostermann \cite{Cassini23} proposed a fundamental acceleration technique for ExpRK integrators and Lawson schemes to efficiently solve semilinear advection-diffusion-reaction equations. The model under consideration is of the form:
\begin{equation}
    \partial_t u(t, x)=\nabla \cdot (a(x) \nabla u(t, x)) + \nabla \cdot (b(x) u(t, x))+r(t, x, u(t, x)),
\end{equation}
where $a(x)$ is the diffusion tensor, $b(x)$ is the velocity field, and $r(t, x, u)$ represents the reaction term. The main computational challenge of exponential integrators lies in efficiently evaluating matrix functions such as $\textrm{e}^{\tau A}$ and $\varphi_1(\tau A)$, whose definition can be obtained from \eqref{phirecursive}.
To improve efficiency, the authors propose choosing a linear operator $A$ with constant coefficients so that the computation of these functions can be performed more efficiently, for example, using FFT in cases with constant diffusion. The operator $A$ is defined as:
\begin{equation}
    Au=\lambda \sum_{\mu=1}^{d} \max_x a_{\mu\mu}(x) \partial^2_{x_\mu}u+ \beta \cdot \nabla u\,,
\end{equation}
where the parameter $\lambda$ is determined to ensure stability based on a spectral analysis of the Laplace operator. For instance, when solving pure diffusion equations, the authors demonstrate that the ETD Euler \eqref{ETDeuler} 
is unconditionally stable for $\lambda \geq 1/2$. Meanwhile, the Lawson-Euler scheme:
\begin{equation}
    u^{n+1}=\textrm{e}^{\tau \lambda \Delta} (u^n + \tau (1-\lambda) \Delta u^n)
\end{equation}
is stable for $\lambda \geq 0.218$, offering greater flexibility in the choice of the operator. This approach is validated numerically on advection-diffusion-reaction equations in 1D, 2D, and 3D, demonstrating significant computational efficiency improvements over IMEX schemes and traditional exponential methods.

Again in 2024, multiple works were published on the subject of order reduction of ExpRK methods, focusing specifically on problems with non-commutative operators. In a first work, Hoang \cite{hoang2024order} studied the application of ExpRK methods to nonlinear parabolic equations with non-commutative operators. The model under consideration is an equation of the form:
\begin{equation}
    \partial_t u - \Delta u = f(\nabla u, u), \quad u(0) = u_0\,,
\end{equation}
where $f$ is a nonlinear function. To understand the numerical behavior of such methods, the analysis focuses on the simplified problem:
\begin{equation}
    u'(t)+Au(t)=Bu(t), \quad u(0)=u_0\,,
\end{equation}
where $A$ is an operator generating an analytic semigroup, while $B$ is \emph{relatively bounded}\footnote{In general, an operator $B$ is said to be relatively bounded with respect to another operator
$A$ if there exist constants $C\ge 0$ and $\alpha\ge0$
\[
\|Bu\|\le C\le \|Au\|+\alpha\|u\|, \quad \forall u\in D(A)\,.
\]
This, in practice, means that $B$ is not necessarily bounded, but its growth is controlled by $A$. This property is crucial in spectral analysis and numerical methods, ensuring stability and well-behaved approximations in the presence of unbounded operators.} with respect to $A$. The method follows an approach where $A$ is treated exactly and $B$ is integrated explicitly. The work focuses on the phenomenon of order reduction, which occurs when the method's convergence is lower than theoretically expected. The analysis shows that, in the presence of non-commutative operators, higher-order methods suffer from accuracy degradation unless suitable regularity conditions on the initial data are met. Through theoretical results and numerical simulations, the study highlights the limitations and possible improvements of exponential methods in handling such equations.

In a second work authored by Dang and Hoang \cite{dang2024avoid}, the authors investigated the issue of order reduction in third-order ExpRK methods applied to problems with non-commutative operators. In particular, they proved that if specific order conditions are not satisfied, the method's accuracy may drop to approximately $2.5$ instead of $3$. To overcome this problem, the authors propose a 4 stages ExpRK method that meets all the necessary conditions to preserve the expected order. The studied method takes the form:
\begin{equation}
\begin{aligned}
    U_{ni}&=\textrm{e}^{-c_i \tau A} u_n + \tau \sum_{j=1}^{i-1} a_{ij}(-\tau A) B U_{nj}\,,\\
    u_{n+1}&=\textrm{e}^{-\tau A}u_n + \tau \sum_{i=1}^{s} b_i(-\tau A)BU_{ni}\,,
\end{aligned}
\end{equation}
with the coefficients $a_{ij}(-\tau A)$ and $b_i(-\tau A)$ defined through linear combinations of the functions $\varphi_k(-\tau A)$. A detailed convergence analysis is conducted, based on a semigroup formulation in a Banach space, showing that adopting higher regularity on the initial data helps to avoid order reduction. Numerical tests confirm the validity of the approach, demonstrating significant improvements in accuracy and stability compared to standard three-stage methods.

In an additional work, Einkemmer, Hoang and Ostermann \cite{einkemmer2024should} studied the efficiency of exponential integrators for advection-dominated problems. They compared Leja and Krylov-based methods for computing matrix exponential actions, evaluating their performance against explicit RK schemes. Their results show that while exponential integrators perform similarly to explicit methods in fully advection-dominated regimes, they offer advantages when parts of the domain exhibit diffusion. Among the tested approaches, Leja-based methods proved more efficient than Krylov iterations, particularly on modern supercomputers where inner product computations are costly. The paper focused on performance modeling: counting core operations to understand theoretical cost ratios between exponential integrators and explicit schemes. The recent complementary work \cite{dang2025memory}, instead, measured real CPU time in numerical experiments. It showed Krylov-based methods can work better for small steps in real CPU usage -- a nuanced result not apparent from the theoretical model.

In 2024, Wang, Zhao, and Liao \cite{wangzhaoliao} presented a unified theoretical framework for analyzing the energy dissipation properties of three classes of RK methods applied to phase-field Crystal (PFC) models, which describe crystalline dynamics at atomic scales. The fundamental model is given by the gradient flow:
\begin{equation}
\partial_t\Phi = M\frac{\delta E}{\delta\Phi},
\end{equation}
where the Swift-Hohenberg free energy functional is:
\begin{equation}
E[\Phi] := \int_{\Omega}\left[\frac{1}{2}\Phi(I + \Delta)^2\Phi + F(\Phi)\right]\textrm{d}x,
\end{equation}
with $F(\Phi) = \frac{1}{4}\Phi^4 - \frac{\varepsilon}{2}\Phi^2$ and $f(\Phi) = -F'(\Phi)$. Two common choices of $M$ lead to different evolution equations: the $L^2$ gradient flow ($M = -I$) results in the Swift-Hohenberg equation, while the $H^{-1}$ gradient flow ($M = \Delta$) leads to the conservative volume PFC equation. A key challenge is ensuring maximum norm boundedness of solutions without assuming global Lipschitz continuity of $f(\Phi)$. To address this, a linear stabilization technique is introduced:
\begin{equation}
L_\kappa\Phi := (I+\Delta)^2\Phi+ \kappa\Phi\,, \quad f_\kappa(\Phi):= \kappa\Phi+f(\Phi)\,,
\end{equation}
where $\kappa \geq 0$ is a stabilization parameter. The study provides a unified framework for analyzing energy dissipation across different RK methods. The main result states that numerical solutions maintain maximum norm boundedness at all stages, provided that the \emph{differentiation matrix} associated with the RK method satisfies a positivity condition. Among the three classes analyzed, ExpRK methods, including explicit and ETD RK methods, are particularly effective, as they treat stiff linear terms exactly while handling nonlinear terms explicitly, improving stability and computational efficiency. The other two classes examined are implicit-explicit Runge-Kutta and Corrected Integrating Factor Runge-Kutta methods. The novelty of the work lies in proving that these methods preserve the original energy dissipation laws if their differentiation matrices are positive definite. This is achieved using a differential form approach, discrete orthogonal convolution kernels, and mathematical induction, paving the way for rigorous energy stability in dissipative semilinear parabolic problems.

\section{Examples}
\label{Sec4}
To conclude our manuscript, we present three examples that illustrate the general functionality of ExpRK methods. To this end, we consider standard methods of two different orders (2 and 4) and compare an ExpRK method with a Rosenbrock method and a classical explicit RK method of the same order. The purpose of these examples is not to draw general conclusions in favor of exponential methods, nor to test specific features or provide an exhaustive assessment of their performance. Rather, our goal is to gain insight into which classes of equation could benefit from the use of ExpRK methods, by means of comparison with well-established classical schemes. In this perspective, the examples are meant to have a primarily illustrative and didactic value, offering a practical introduction to the behavior of ExpRK methods in representative settings.

We remark that, in 2020, Montanelli and Bootland \cite{montanelli2020solving} published a comprehensive review of exponential methods, comparing various formulas for solving stiff PDEs in multiple dimensions. Their analysis, focused on periodic boundary conditions and semilinear stiff PDEs, concluded that the ETDRK4 scheme by Cox and Matthews \cite{cox2002exponential}  (denoted as ExpRK4 in the following, for ease and uniformity of notation) remains one of the most effective approaches despite the proliferation of more complex alternatives. Their work is particularly notable for its examination of both diffusive problems (characterized by large negative eigenvalues) and dispersive problems (distinguished by their highly oscillatory nature), as well as its discussion of computational techniques for the $\varphi$-functions based on contour integrals. In line with Montanelli and Bootland’s review \cite{montanelli2020solving}, we restrict our attention to explicit exponential methods. These methods are by far the most widely employed in practice: they often yield superior performance when stiffness is predominantly carried by the linear dynamics while the nonlinear contribution remains comparatively mild. Moreover, they avoid the extra complications associated with Newton-type iterations, such as the need to solve nonlinear systems at every step and to tune iterative solver tolerances; these steps bring additional computational cost and potential convergence issues. For these reasons, explicit exponential formulations are a natural and pragmatic choice in many applications. Nevertheless, in this work we also deliberately test them on less typical, small-scale examples to probe their behavior in unconventional settings.
\smallskip

The three didactic examples we propose have dimension 1, 2 and 50 (where here the `dimension' refers to the size of the system of ODEs), respectively. Specifically, for the 1D case, we revisit the introductory example presented by Cox and Matthews in 2002 \cite[Sect. 4.1]{cox2002exponential}; for the 2D case, we examine a classical model from nonlinear mechanics and dynamical systems, namely the Duffing oscillator\footnote{For a comprehensive analysis of the nonlinear Duffing oscillator and its applications, we recommend the following volume \cite{DuffingBook} to the interested reader.}; the last example is the result of applying the pseudospectral discretization method \cite{bdgsv16} to an Ikeda DDE \cite{ELLG04}. In what follows, we briefly introduce three order $2$ methods and three order $4$ methods for solving stiff systems of ODEs. As already noticed, stiff ODEs require specialized numerical techniques due to their rapidly changing components existing alongside slowly evolving ones. The six methods under consideration %--RK2, ExpRK2, Rb2, RK4, ExpRK4 and Rb4--
represent different approaches to handling stiffness, while all maintaining second-order [resp. fourth-order] accuracy. Their properties suggest different performance characteristics when applied to stiff problems.

For each test case examined, we conduct a convergence analysis by evaluating the relative error against the analytical solution at the final time for various time steps. We recall that the relative error is calculated as:

\begin{equation}
\text{Error} = \left|\frac{u_{\text{numerical}} - u}{u}\right|.
\end{equation}

In the first example we present, there exists an analytical expression for the solution $u$ of the one-dimensional ODE, and we use this exact expression in the formula above. In the second two-dimensional example, such an explicit, exact expression does not exist, and we compare the methods against a reference expression for the solution $u$. This is also the case for the third example, for which we compute the reference solution using a built-in MATLAB integrator for DDEs. Additionally, we measure the computational time required by each method, enabling a comprehensive comparison of both accuracy and efficiency. This dual assessment, as emphasized earlier, is particularly important when dealing with stiff systems, where the choice of numerical method significantly impacts the simulation's precision and performance. Through this analysis, we aim to further illustrate the practical considerations discussed in previous sections and provide additional insights into the behavior of these methods across different problem types and stiffness regimes. All simulations were performed in MATLAB, version R2024a, on a 64-bit Microsoft Windows system.
\medskip

\paragraph{\sc The Six Methods}

Here, we describe the six numerical methods which we have used for our simulations, namely for solving (in one case non-autonomous, in the other cases autonomous) systems of ODEs of the form:
\begin{equation}
\frac{d\mathbf{u}}{dt} = \mathbf{F}(t, \mathbf{u})\,.
\end{equation}
In particular, we consider three methods of order 2 and three methods of order 4.

\subsection*{Explicit Runge-Kutta 2 Method (RK2)}

RK2 updates the solution as follows:
\begin{align*}
\mathbf{k}_1 &= \mathbf{F}(t_n, \mathbf{u}_n), \\
\mathbf{k}_2 &= \mathbf{F}\left(t_n + \frac{h}{2}, \mathbf{u}_n + \frac{h}{2} \mathbf{k}_1\right), \\
\mathbf{u}_{n+1} &= \mathbf{u}_n + h \mathbf{k}_2.
\end{align*}

\subsection*{Exponential Runge-Kutta 2 Method (ExpRK2)}
We can rewrite the equation as follows:
\begin{equation}
\frac{d\mathbf{u}}{dt} = A\mathbf{u} + \mathbf{g}(t, \mathbf{u})\,,
\end{equation}
where $A\mathbf{u}$ is the linear term. The ExpRK2 method updates the solution as follows:
\begin{align*}
\mathbf{a} &=\e^{h A} \mathbf{u}_n + h \varphi_1(h A) \mathbf{g}(t_n, \mathbf{u}_n)\,\\
\mathbf{u}_{n+1} &= \mathbf{a} + h \varphi_2(h A) (\mathbf{g}(t_n+h, \mathbf{a})-\mathbf{g}(t_n, \mathbf{u}_n))\,,
\end{align*}
where $\varphi_1$ and $\varphi_2$ are defined from \eqref{phirecursive}.

\subsection*{Rosenbrock 2 Method (Rb2)}
The Rosenbrock method is an implicit method designed for stiff ODEs, incorporating a Jacobian-based correction. The update rule is given by:
\begin{align*}
\mathbf{k} &= \left( I - \gamma h \mathbf{J}_n \right)^{-1} \mathbf{F}(t_n, \mathbf{u}_n)\,, \\
\mathbf{u}_{n+1} &= \mathbf{u}_n + h \mathbf{k}\,,
\end{align*}
where  
\[
\mathbf{J}_n=\frac{\partial \mathbf{F}}{\partial \mathbf{u}} \Bigg|_{(t_n,\mathbf{u}_n)}
\]
is the Jacobian matrix, and $\gamma$ is a stabilization parameter; in particular, $\gamma=0.5$ provides a method of order 2.

\subsection*{Explicit Runge-Kutta 4 Method (RK4)}
The classical order 4 method that we consider is RK4, which we already introduced in Section \ref{Sec22}.

\subsection*{Exponential Runge-Kutta 4 Method (ExpRK4)}
Among the order 4 exponential methods, we consider ExpRK4 \cite[eq. (26)-(29)]{cox2002exponential}. Considering again the equation,
\begin{equation}
\frac{d\mathbf{u}}{dt} = A\mathbf{u} + \mathbf{g}(t, \mathbf{u})\,,
\end{equation}
ExpRK4 updates the solution as follows:
\begin{align*}
\mathbf{g}_1&=\mathbf{g}(t_n, \mathbf{u}_n),\,\\
\mathbf{a} &=\e^{h A/2} \mathbf{u}_n + \frac{h}{2} \varphi_1(h A/2)\mathbf{g}_1, \,\\
\mathbf{g}_2&=\mathbf{g(t_n+h/2,\mathbf{a}}),\,\\
\mathbf{b} &=\e^{hA/2}\mathbf{u}_n +\frac{h}{2} \varphi_1(h A/2)\mathbf{g}_2,\,\\
\mathbf{g}_3&=\mathbf{g(t_n+h/2,\mathbf{b}}),\,\\
\mathbf{c} &=\e^{hA/2}\mathbf{a} +\frac{h}{2}(2\mathbf{g}_3-\mathbf{g}_1),\,\\
\mathbf{g}_4&=\mathbf{g}(t_n+h, \mathbf{c})),\,\\
b_1(hA)&=\varphi_1(hA)-3\varphi_2(hA)+4\varphi_3(hA),\,\\
b_2(hA)&=2\varphi_2(h A)-4\varphi_3(hA),\,\\
b_3(hA)&=-\varphi_2(hA)+4\varphi_3(hA),\,\\
\mathbf{u}_{n+1} &= \e^{hA}\mathbf{u}_n + h (b_1(hA)\mathbf{g}_1+b_2(hA) (\mathbf{g}_2+\mathbf{g}_3)+b_3(hA)\mathbf{g}_4),\,
\end{align*}
where $\varphi_1$, $\varphi_2$ and $\varphi_3$ are defined from \eqref{phirecursive}.

\subsection*{Rosenbrock 4 Method (Rb4)}
The order 4 Rosenbrock-Wanner method we consider is defined by four stages. Introducing for ease of notation $f(t_n,\mathbf{u}_n):=f(\mathbf{u}_n)$, the update rule is given by:
\begin{align*}
    &\left(I + \gamma_{ii}h\mathbf{J}_n\left(\mathbf{u}_n\right)\right)\mathbf{k}_i \\
    &= hf\left(\mathbf{u}_n +\sum_{j=1}^{i-1} \alpha_{ij}\mathbf{k}_j \right) + h\mathbf{J}_n(\mathbf{u}_n)\sum_{j=1}^{i-1}\gamma_{ij}\mathbf{k}_j,\quad i=1,\ldots 4,\nonumber\\
    \mathbf{u}_{n+1}&=\mathbf{u}_n+\sum_{i=1}^{4} \beta_i\mathbf{k}_i,
\end{align*}
for constant $\alpha_{ij},\,\beta_i,\,\gamma_{ij}$ which are chosen in order to satisfy the order 4 conditions (see, e.g., \cite[Section IV.7]{hw91}).
\smallskip

\subsection*{\sc 1D. C\&M Starting Example}
Following the same approach presented by Cox and Matthews \cite{cox2002exponential} in their analysis of ETD methods, we consider a non-autonomous ODE model with rapid decay, characteristic of stiff systems:
\begin{equation}
u'(t) = ku(t) + \sin(t), \quad u(0) = u_0 =1 \,,
\end{equation}
where we have chosen $k = -100$ to generate the rapid linear decay typical of stiff problems, in contrast with the $\mathcal{O}(1)$ time scale of the forcing term. This equation, while simple in form, captures the essential features of stiff systems: a fast decaying homogeneous solution combined with a (relatively) slow varying forcing term. For subsequent evaluation of the different numerical schemes, we utilize the exact analytical solution given by:
\begin{equation}
u(t) = -\frac{-e^{kt}(2 + k^2) + \cos(t) + k\sin(t)}{1 + k^2}.
\end{equation}
The behavior of this differential equation can be divided into two distinct phases: an initial fast transient phase, during which the solution rapidly approaches the slow manifold, followed by a phase in which the solution evolves along this slow manifold. As noticed in \cite{cox2002exponential}, when $k < 0$ and $|k| \gg 1$, the solution quickly approaches a slow manifold on which
\[
u \sim -\frac{\sin(t)}{k} - \frac{\cos(t)}{k^2} + \ldots\,.
\]
A robust numerical method must accurately capture both of these phases and work effectively when the time step $h$ is of the order of $1/|k|$. This test case allows us to evaluate how well our three order $2$ and three order $4$ methods handle the challenges posed by stiffness, particularly the balance between stability and accuracy when using relatively large time steps compared to the fastest time scale in the system. We examined the performance of the methods over the interval $[0, \pi/2]$ using various time steps ranging from $h = 0.0001$ to $h = 0.1$, comparing both the accuracy of the numerical solutions against the analytical solution and the computational efficiency of each method. To estimate the elapsed time, we performed each simulation several times (1000 in the first test, 10 in the second and 100 in the third) and reported the average.
\bigskip

\begin{table}[ht!]
\centering
\caption{Numerical methods performance for various step sizes (order 2 methods).}
\label{Table:ex1}
\resizebox{\linewidth}{!}{\begin{tabular}{lcccccc}
\toprule
$h$ & RK2 Error & ExpRK2 Error & Rb2 Error & RK2 Time (s) & ExpRK2 Time (s) & Rb2 Time (s) \\
\midrule
$1 \times 10^{-4}$ & $1.2551 \times 10^{-9}$ & $8.3276 \times 10^{-10}$ & $2.3212 \times 10^{-9}$ & 0.05355 & 0.05829 & 0.05625 \\
$5 \times 10^{-4}$ & $3.1432 \times 10^{-8}$ & $2.0438 \times 10^{-8}$ & $3.1446 \times 10^{-8}$ & 0.01183 & 0.01188 & 0.01070 \\
$1 \times 10^{-3}$ & $1.2779 \times 10^{-7}$ & $8.0986 \times 10^{-8}$ & $1.6274 \times 10^{-7}$ & 0.00556 & 0.00584 & 0.00585 \\
$5 \times 10^{-3}$ & $3.8530 \times 10^{-6}$ & $1.9199 \times 10^{-6}$ & $4.4691 \times 10^{-6}$ & 0.00118 & 0.00131 & 0.00118 \\
$1 \times 10^{-2}$ & $2.3082 \times 10^{-5}$ & $7.5741 \times 10^{-6}$ & $2.1139 \times 10^{-5}$ & 0.00059 & 0.00065 & 0.00061 \\
$5 \times 10^{-2}$ & $7.0246 \times 10^{30}$ & $4.8629 \times 10^{-5}$ & $1.4075 \times 10^{-4}$ & 0.00010 & 0.00013 & 0.00011 \\
$1 \times 10^{-1}$ & $2.9523 \times 10^{27}$ & $2.5459 \times 10^{-4}$ & $4.1347 \times 10^{-2}$ & 0.00005 & 0.00006 & 0.00005 \\
\bottomrule
\end{tabular}}%}

\end{table}

\begin{table}[ht!]
\centering
\caption{Numerical methods performance for various step sizes (order 4 methods).}
\label{Table:ex1_order4}
\resizebox{\linewidth}{!}{\begin{tabular}{lcccccc}
\toprule
$h$ & RK4 Error & ExpRK4 Error & Rb4 Error & RK4 Time (s) & ExpRK4 Time (s) & Rb4 Time (s) \\
\midrule
$1 \times 10^{-4}$ & $1.0583 \times 10^{-14}$ & $4.4262 \times 10^{-10}$ & $7.3560 \times 10^{-14}$ & 0.05772 & 0.07283 & 0.17238 \\
$5 \times 10^{-4}$ & $6.5026 \times 10^{-12}$ & $1.9799 \times 10^{-12}$ & $4.2579 \times 10^{-11}$ & 0.01216 & 0.01484 & 0.03288 \\
$1 \times 10^{-3}$ & $1.0440 \times 10^{-10}$ & $8.0308 \times 10^{-13}$ & $6.2264 \times 10^{-10}$ & 0.00576 & 0.00692 & 0.01561 \\
$5 \times 10^{-3}$ & $7.6433 \times 10^{-8}$ & $6.7939 \times 10^{-13}$ & $8.0284 \times 10^{-6}$ & 0.00130 & 0.00150 & 0.00327 \\
$1 \times 10^{-2}$ & $1.5387 \times 10^{-6}$ & $1.1421 \times 10^{-11}$ & $5.7421 \times 10^{-6}$ & 0.00061 & 0.00073 & 0.00161 \\
$5 \times 10^{-2}$ & $6.4073 \times 10^{36}$ & $3.0104 \times 10^{-9}$ & $1.7416 \times 10^{-4}$ & 0.00011 & 0.00013 & 0.00029 \\
$1 \times 10^{-1}$ & $5.8631 \times 10^{40}$ & $4.7821 \times 10^{-8}$ & $5.2099 \times 10^{-4}$ & 0.00006 & 0.00007 & 0.00016 \\
\bottomrule
\end{tabular}}%}
\end{table}

\begin{figure}[ht!]
\centering
\begin{subfigure}{.475\textwidth}
  \centering
  \includegraphics[width=1\linewidth]{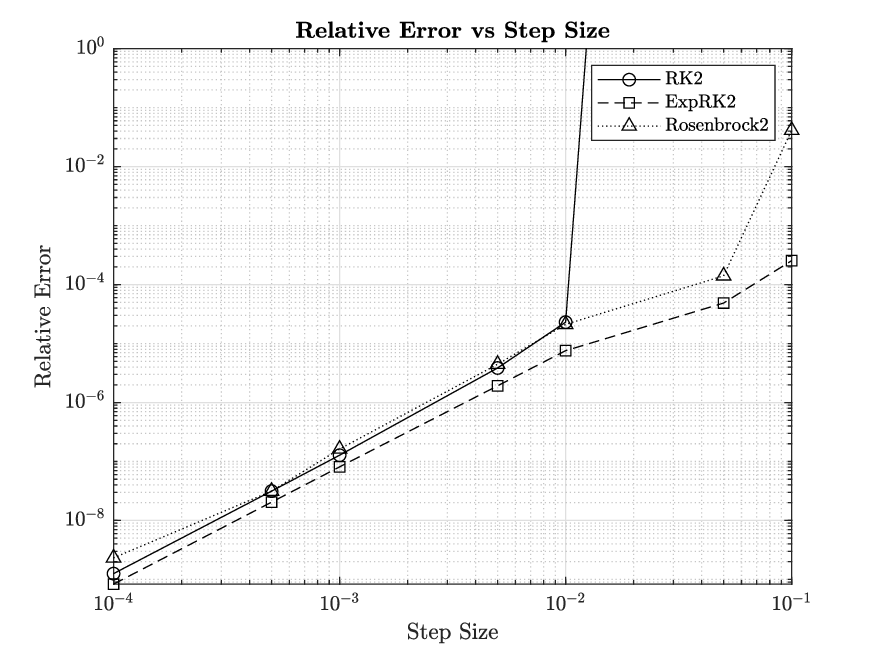}
    \caption{}
  \label{fig:errortime_plot1}
\end{subfigure}%
\begin{subfigure}{.475\textwidth}
  \centering
  \includegraphics[width=1\linewidth]{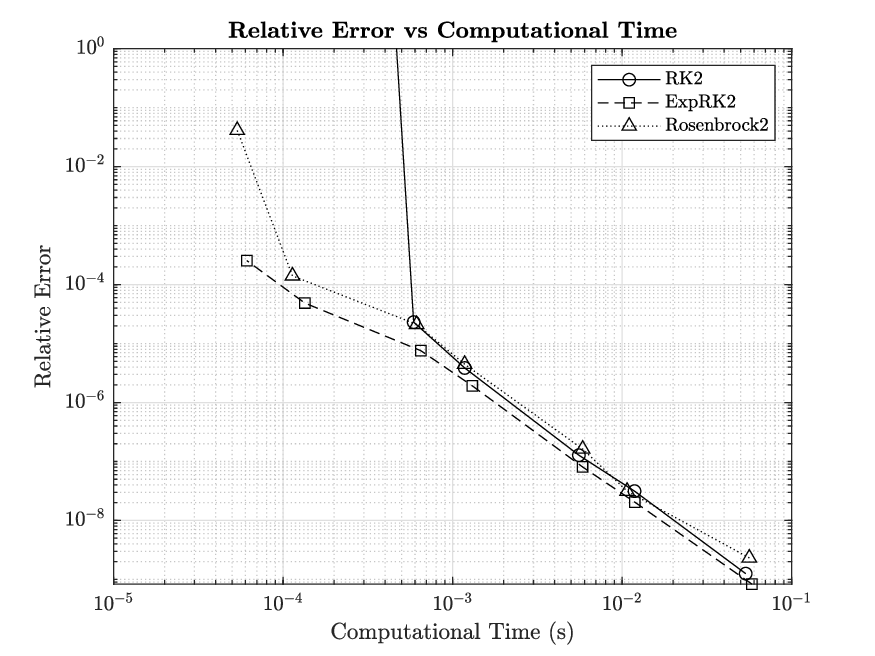}
  \caption{}
  \label{fig:errortime_plot2}
\end{subfigure}\\

\begin{subfigure}{.475\textwidth}
  \centering
  \includegraphics[width=1\linewidth]{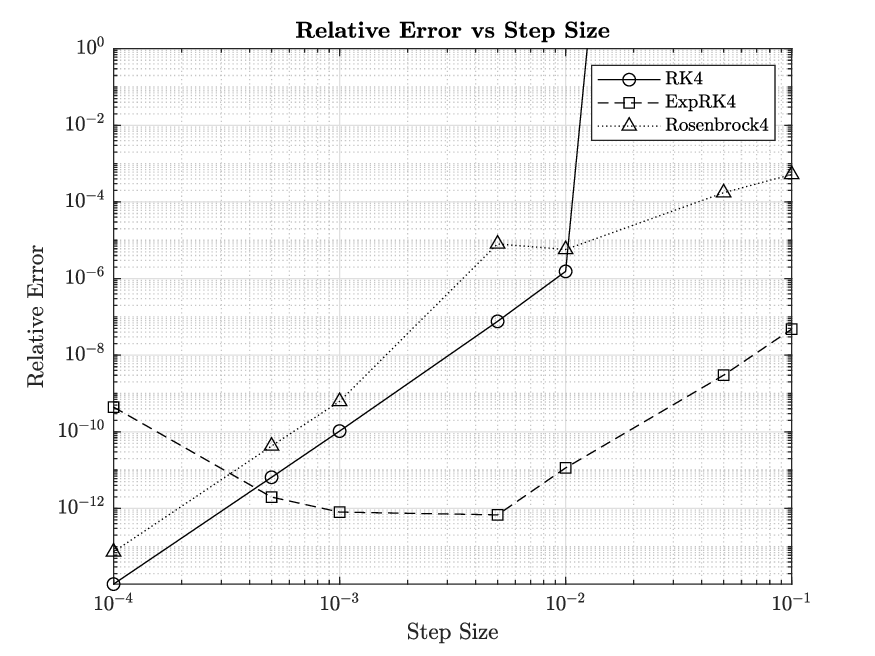}
    \caption{}
  \label{fig:errortime_plot1_order4}
\end{subfigure}%
\begin{subfigure}{.475\textwidth}
  \centering
  \includegraphics[width=1\linewidth]{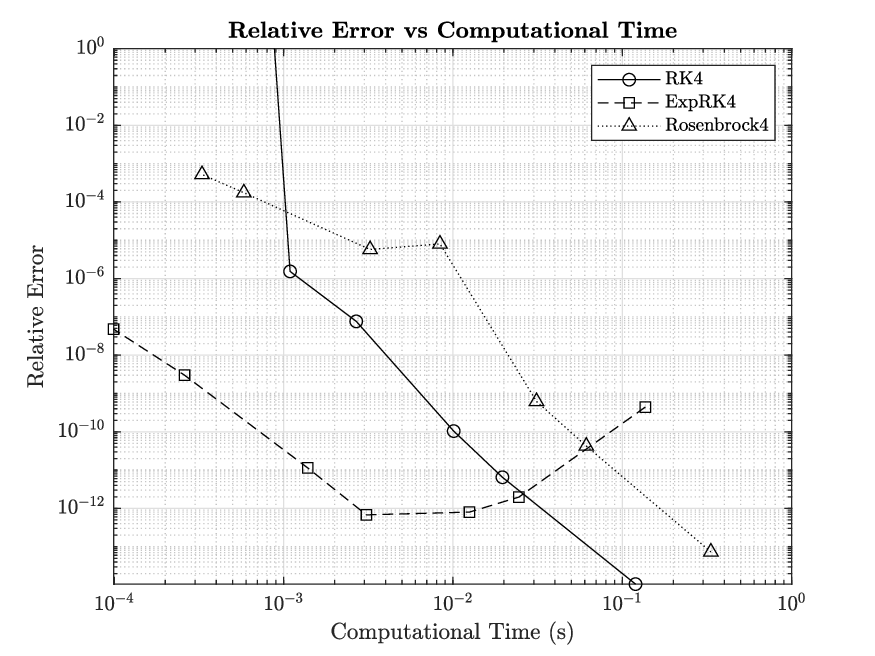}
  \caption{}
  \label{fig:errortime_plot2_order4}
\end{subfigure}
\caption{Relative error as a function of step size, showing how accuracy of each method of order 2
({\subref{fig:errortime_plot1}}) and order 4 (\subref{fig:errortime_plot1_order4}) evolves with different step choices and relative error as a function of computational time, highlighting the efficiency  of methods of order 2 ({\subref{fig:errortime_plot2}}) and order 4 ({\subref{fig:errortime_plot2_order4}})}
\label{fig:errortime_plot}
\end{figure}

The results presented in Figure \ref{fig:errortime_plot} and Tables \ref{Table:ex1}-\ref{Table:ex1_order4} highlight the distinct behaviors of the three order $2$ and the three order $4$ (respectively) numerical methods under analysis. 

RK2 and RK4 demonstrate stability only for small step sizes, specifically for values of $h \leq 0.01$. However, for $h \geq 0.05$, the relative error grows uncontrollably. This behavior confirms a strict stability limit around $h = 0.01$, beyond which the method fails to provide reliable results. On the other hand, both exponential methods and Rb4 exhibit excellent stability across all tested step sizes. Even for larger ones, the methods maintain a reasonably small relative error. Rb2 also demonstrates stable behavior, albeit within a more restricted range. For step sizes up to $0.01$, the method produces accurate results. However, for $h = 0.05$, the relative error starts increasing significantly, and for $h = 0.1$ it becomes quite large ($1.53 \times 10^{-1}$). The superior accuracy makes ExpRK2 a strong candidate for solving stiff problems among those of order 2. With respect to those of order 4, ExpRK4 provides much more accurate results than its Rosenbrock counterpart for most values of $h$, although Rb4 behaves well for all tested step sizes.

Regarding computational efficiency, differences in computation times among the methods of order 2 are negligible, making observations about performance insignificant. In contrast, ExpRK4 appears to be the most efficient among the methods of order 4. It is worth noting that a relative error close to machine precision is already reached using $h=0.005$, which makes the use of smaller step sizes not sensible for this example.
\smallskip

\paragraph{\bf Final considerations}

Considering both accuracy and efficiency, ExpRK2 emerges as the best-balanced method among those of order 2. It sustains high accuracy even with large step sizes, allowing for significantly larger steps compared to the other methods. Although it has a slight computational overhead per step, the reduced number of required steps makes it more efficient overall. 

Rb2 offers similar computation times for very small step sizes but with stricter stability constraints: the method appears to work very effectively for a problem of small scale. RK2, in contrast, proves to be the least suitable method for this stiff problem, as it requires extremely small time steps to maintain stability and lacks competitiveness in both accuracy and computational efficiency.

These findings align with the conclusions drawn by Cox and Matthews in \cite{cox2002exponential}, confirming the superiority of exponential methods for solving stiff problems. In particular, ExpRK2 offers a significant advantage in terms of stability and accuracy, permitting the use of much larger step sizes than conventional methods. The superiority of ExpRK4 with respect to the tested methods of the same order is even more evident. For stiff problems such as the one studied here, characterized by a parameter $k=-100$, exponential methods represent the optimal choice among the three methods of the corresponding order, striking an ideal balance between accuracy, stability, and computational efficiency.

\subsection*{\sc 2D. Duffing Oscillator} The system under consideration in this second example is the Duffing oscillator, a fundamental model in nonlinear dynamics describing a forced harmonic oscillator with a cubic nonlinearity in the restoring force. The governing differential equation is given by:
\begin{equation}
    u'' + \omega u + k u^3 = 0\,,
\end{equation}
where $u$ represents the position, $\omega$ is the natural frequency, and $k$ is the nonlinear coefficient. To analyze its temporal evolution, we rewrite the system as a first-order autonomous system by introducing the state variable $v = u'$, leading to the phase-space formulation:

\begin{equation}
    \begin{cases}
        u' = v\,, \\
        v' = -\omega u - k u^3\,.
    \end{cases}
\end{equation}
This formulation allows us to study the system as a two-dimensional dynamical problem. For numerical integration, it is convenient to express this system in a matrix-vector form:
\begin{equation}
    \begin{bmatrix}
        u' \\
        v'
    \end{bmatrix}
    = 
    \begin{bmatrix}
        0 & 1 \\
        -\omega & 0
    \end{bmatrix}
    \begin{bmatrix}
        u \\
        v
    \end{bmatrix}
    +
    \begin{bmatrix}
        0 \\
        -k u^3
    \end{bmatrix}\,.
\end{equation}
This matrix formulation highlights the linear and nonlinear contributions to the system's dynamics. The stiffness of the system becomes significant for large values of $k$, where the nonlinear term induces rapid variations in the solutions. This necessitates the use of numerical schemes capable of handling both slow dynamics and the sharp transitions associated with high-curvature regions in the phase space trajectories.
\bigskip

\begin{table}[ht!]
\centering
\caption{Numerical methods performance for various step sizes (order 2 methods).}
\label{Table:duffing_o2}
\resizebox{\linewidth}{!}{\begin{tabular}{lcccccc}
\toprule
$h$ & RK2 Error & ExpRK2 Error & Rb2 Error & RK2 Time (s) & ExpRK2 Time (s) & Rb2 Time (s) \\
\midrule
$1 \times 10^{-4}$ & $6.2072 \times 10^{-5}$ & $6.0105 \times 10^{-6}$ & $8.7547 \times 10^{-5}$ & 0.09026 & 0.54986 & 1.13880 \\
$5 \times 10^{-4}$ & $1.9363 \times 10^{-3}$ & $2.5767 \times 10^{-4}$ & $2.2968 \times 10^{-3}$ & 0.01701 & 0.09588 & 0.20080 \\
$1 \times 10^{-3}$ & $9.6628 \times 10^{-3}$ & $1.5697 \times 10^{-3}$ & $9.7327 \times 10^{-3}$ & 0.00866 & 0.04709 & 0.10134 \\
$5 \times 10^{-3}$ & $6.3873 \times 10^{-1}$ & $1.4548 \times 10^{-1}$ & $3.4409 \times 10^{-1}$ & 0.00176 & 0.00874 & 0.01715 \\
$1 \times 10^{-2}$ & $2.9671 \times 10^{0}$ & $9.5176 \times 10^{-1}$ & $1.1835 \times 10^{0}$ & 0.00092 & 0.00674 & 0.01204 \\
$5 \times 10^{-2}$ & NaN                  & $1.3049 \times 10^{0}$ & $1.2752 \times 10^{-1}$ & 0.00008 & 0.00224 & 0.00277 \\
$1 \times 10^{-1}$ & NaN                  & $1.6391 \times 10^{0}$ & $8.3403 \times 10^{-1}$ & 0.00002 & 0.00059 & 0.00113 \\
\bottomrule
\end{tabular}}
\end{table}

\begin{table}[ht!]
\centering
\caption{Numerical methods performance for various step sizes (order 4 methods).}
\label{Table:duffing_o4}
\resizebox{\linewidth}{!}{\begin{tabular}{lcccccc}
\toprule
$h$ & RK4 Error & ExpRK4 Error & Rb4 Error & RK4 Time (s) & ExpRK4 Time (s) & Rb4 Time (s) \\
\midrule
$1 \times 10^{-4}$ & $4.7521 \times 10^{-12}$ & $1.0584 \times 10^{-9}$ & $6.0232 \times 10^{-11}$ & 0.83861 & 15.54832 & 2.15542 \\
$5 \times 10^{-4}$ & $6.7006 \times 10^{-9}$  & $2.3246 \times 10^{-9}$ & $2.2207 \times 10^{-8}$  & 0.16769 & 0.80160  & 0.50081 \\
$1 \times 10^{-3}$ & $1.3075 \times 10^{-7}$  & $3.8983 \times 10^{-8}$ & $2.0253 \times 10^{-7}$  & 0.06966 & 0.20216  & 0.20446 \\
$5 \times 10^{-3}$ & $1.9809 \times 10^{-4}$  & $4.9790 \times 10^{-5}$ & $6.3571 \times 10^{-4}$  & 0.01447 & 0.02514  & 0.04033 \\
$1 \times 10^{-2}$ & $5.4686 \times 10^{-3}$  & $1.3051 \times 10^{-3}$ & $2.5171 \times 10^{-2}$  & 0.00977 & 0.01663  & 0.02508 \\
$5 \times 10^{-2}$ & $4.1232 \times 10^{-1}$  & $3.2182 \times 10^{0}$  & $4.6251 \times 10^{-1}$  & 0.00217 & 0.00512  & 0.00734 \\
$1 \times 10^{-1}$ & $1.9771 \times 10^{0}$   & NaN                     & $3.2802 \times 10^{-2}$ & 0.00105 & 0.00300  & 0.00378 \\
\bottomrule
\end{tabular}}
\end{table}

\begin{figure}[ht!]
\centering
\begin{subfigure}{.475\textwidth}
  \centering
  \includegraphics[width=1\linewidth]{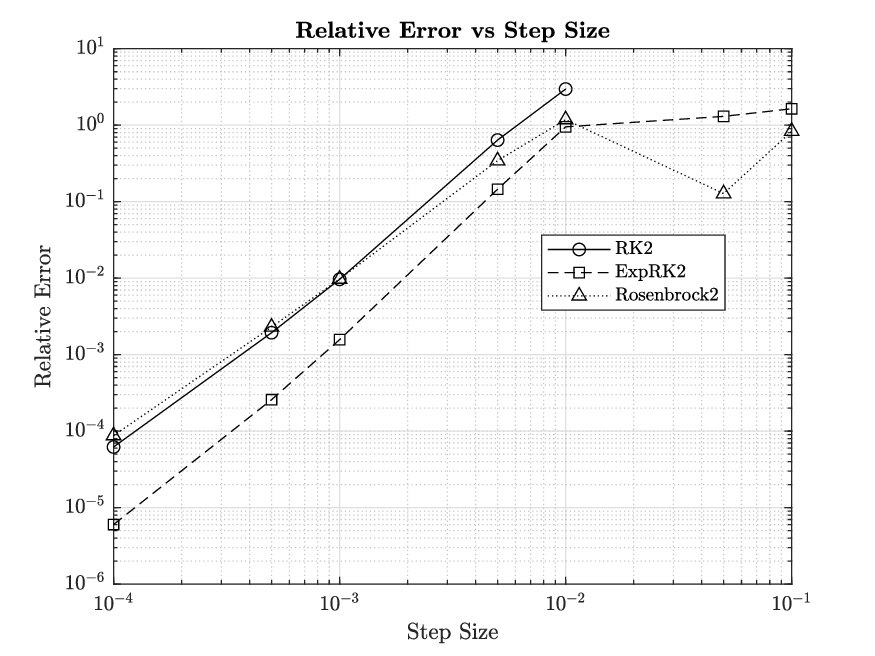}
  \caption{}
  \label{fig:error_step_plot}
\end{subfigure}%
\begin{subfigure}{.475\textwidth}
  \centering
  \includegraphics[width=1\linewidth]{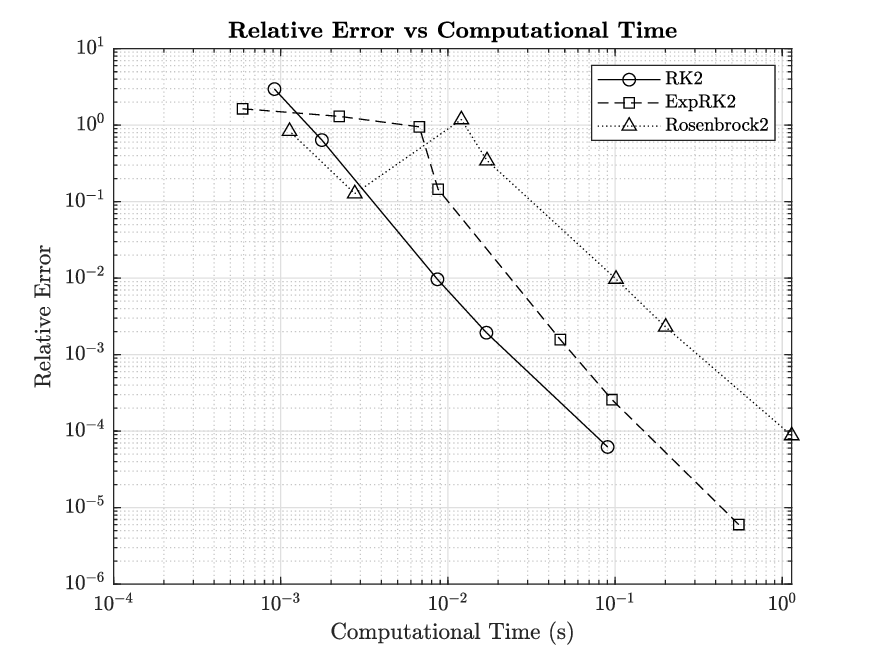}
  \caption{}
  \label{fig:error_time_plot}
\end{subfigure}
\begin{subfigure}{.475\textwidth}
  \centering
  \includegraphics[width=1\linewidth]{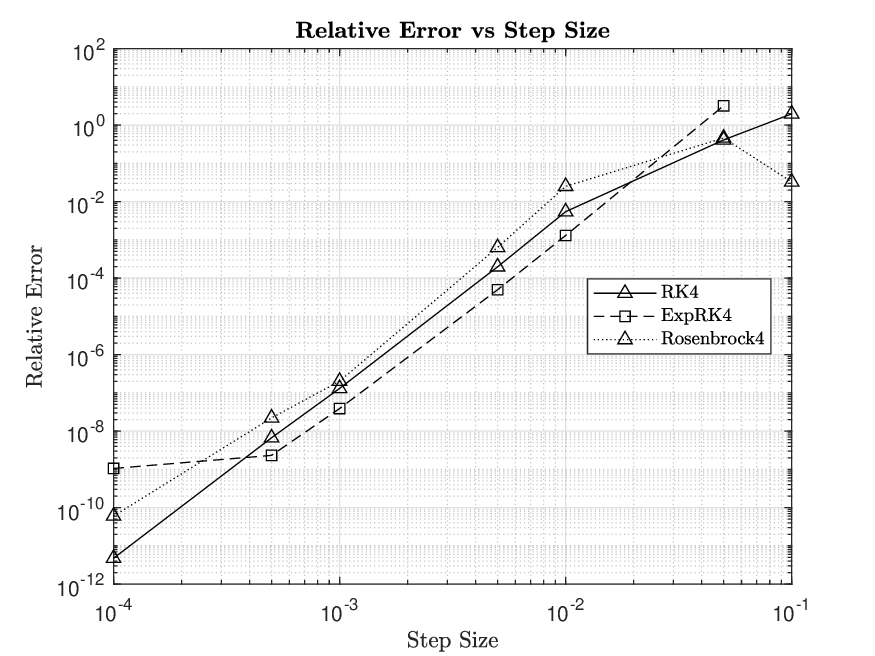}
    \caption{}
  \label{fig:errortime_plot2D_1_order4}
\end{subfigure}%
\begin{subfigure}{.475\textwidth}
  \centering
  \includegraphics[width=1\linewidth]{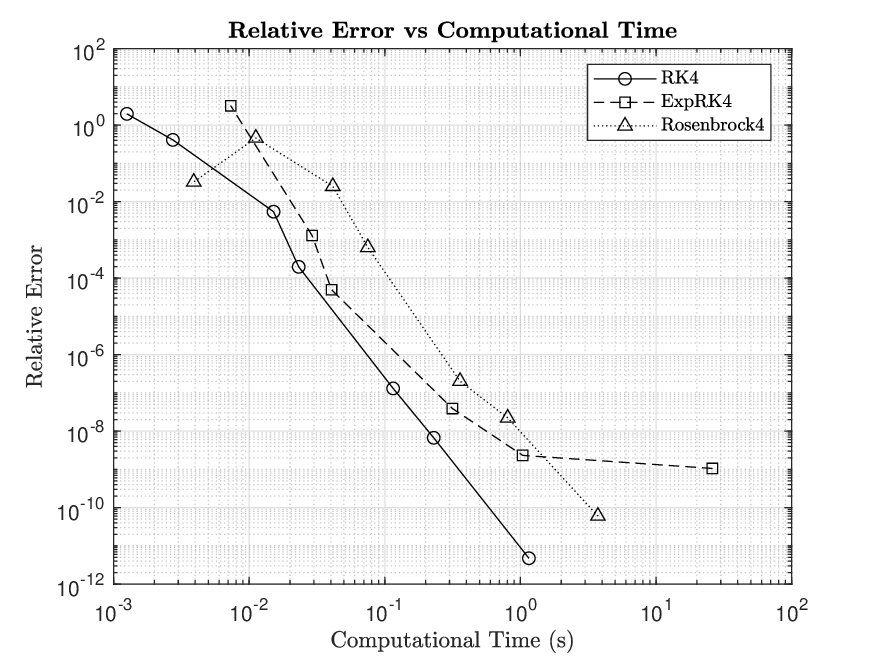}
  \caption{}
\label{fig:errortime_plot2D_2_order4}
\end{subfigure}
\caption{Relative error as a function of step size, showing how accuracy of each method of order 2 ({\subref{fig:error_step_plot}}) and order 4 ({\subref{fig:errortime_plot2D_1_order4}}) evolves with different step choices and relative error as a function of computational time, highlighting the efficiency  of methods of order 2 ({\subref{fig:error_time_plot}}) and order 4 ({\subref{fig:errortime_plot2D_2_order4}})}
\label{fig:duffing_plots}
\end{figure}

The results presented in Figure \ref{fig:duffing_plots} and Tables \ref{Table:duffing_o2}-\ref{Table:duffing_o4} provide a comprehensive analysis of six numerical methods of order 2 and 4 applied to the Duffing oscillator system\footnote{All relative errors are computed with respect to a high-precision numerical solution obtained using MATLAB's \texttt{ode15s} solver, as no analytical solution is available for this nonlinear system.} integrated over the time interval $[0,20]$ with $\omega = 1$ and $k = 100$.

RK2 demonstrates severe stability limitations when applied to this problem. For step sizes $h \leq 0.001$, the method produces reasonable errors, but as $h$ increases to $0.01$, the relative error grows to approximately $2.97$ and, most importantly, for larger step sizes ($h \geq 0.05$), the method fails entirely, producing \texttt{NaN} values. On the other hand, stability constraints on RK4 are not as evident. This suggests that the problem, with these parameters, can be considered only midly stiff.

Rb2 exhibits significantly better stability characteristics than RK2: it maintains reasonable accuracy for step sizes up to $h = 0.01$, with errors remaining on the order of $10^0$. 

ExpRK2 demonstrates superior accuracy for most step sizes up to $h = 0.01$, particularly excelling at moderate step sizes where it maintains errors significantly lower than the other methods. For $h = 5 \times 10^{-4}$ and $h = 1 \times 10^{-3}$, ExpRK2 shows errors approximately an order of magnitude smaller than both RK2 and Rb2, despite not being designed so as to capture stiffness phenomena caused by the \emph{nonlinear} terms. However, for larger step sizes ($h \geq 0.05$), the errors increase substantially, exceeding the error of Rb2 at these larger step sizes. ExpRK4 exhibits similar behavior, in that it provides more accurate results for the intermediate step sizes. However, the difference in this case is not as significant.

Regarding computational efficiency, across all step sizes small enough for RK2 to obtain a solution, RK2 is the fastest method. ExpRK2 occupies an intermediate position, with computation times around 3 times higher than RK2 for a given (relative) error. Rb2 emerges as the most computationally intensive method, requiring instead around 10 times the time required by RK2 for a given error. Among the methods of order 4, the classical RK is the fastest again, but the exponential method ExpRK4 is faster than the Rosenbrock method of the same order for intermediate step sizes. In any case, the difference among the three methods in terms of performance is not as significant as within the methods of order 2.

\smallskip

\paragraph{\bf Final considerations}

The analysis of these numerical methods for the Duffing oscillator reveals important trade-offs between accuracy, stability, and computational efficiency. RK2, while computationally inexpensive, suffers from severe stability limitations that render it unsuitable for step sizes larger than $0.001$. Rb2 is the most accurate for step sizes larger than $0.01$, while it becomes less accurate than ExpRK2 for smaller step sizes. Its non-monotonic error behavior at larger step sizes warrants further investigation but could potentially be exploited in specific applications.
As for the methods of order 4, the Rosenbrock one appears again to be the most accurate for larger step sizes, while ExpRK4 may be a better choice for applications with higher accuracy constraints.

ExpRK2 demonstrates excellent accuracy, particularly at larger step sizes where it significantly outperforms both RK2 and Rb2. Indeed, it might be the best combination in terms of accuracy and efficiency, depending on the specific application. The same can be inferred about ExpRK, when compared to the other methods of order 4.

For the Duffing oscillator with the given parameters, the exponential methods appear to offer the most favorable combination of stability, accuracy, and computational efficiency among the methods of the same order. They maintain acceptable accuracy even at large step sizes and does so with reasonable computational cost. This makes them particularly suitable for applications where a balance between accuracy and efficiency is required. However, in applications where low latency is of outmost importance, the Rosenbrock methods may be preferred for moderate step sizes where its superior accuracy is most evident.

These findings highlight the importance of considering the specific characteristics of the differential system when selecting a numerical method, as the relative performance of different methods can vary significantly depending on the nature and parameters of the problem being solved.

\subsection*{\sc Ikeda DDE} 

Lastly, we consider the Ikeda DDE, which models a variety of passive optical systems subject to a delayed feedback \cite{ELLG04}. The general form of the DDE reads
\begin{equation*}
x'(t)=-\lambda x(t)+f(\mu,x(t-\tau)),
\end{equation*}
with $\lambda$ and $\mu$ positive parameters, and $f$ a continuous function. We analyze the following Cauchy problem, obtained by choosing $f(\mu,y)=\mu(1-\text{sin}(y))$ and $\tau=\frac{\pi}{2}$:
\begin{equation}\label{ikeda}
\left\{
\begin{aligned}
&x'(t)=-\lambda x(t)+\mu\left(1-\text{sin}\left(x\left(t-\frac{\pi}{2}\right)\right)\right), \\
&x(t)=\cos(t),\quad t\in\left[-\frac{\pi}{2},0\right].
\end{aligned}
\right.
\end{equation}
In order to make \eqref{ikeda} amenable to our analysis, we need to first discretize the DDE into a finite number $M$ of ODEs. We opt in particular for the \emph{pseudospectral discretization} technique \cite{bdgsv16}. In fact, the ODEs obtained using this approach are defined by a semilinear right-hand side with stiff linear part, regardless of the right-hand side of the DDE. The resulting non-linear term is defined by only one nonzero component, which corresponds to the non-linear term in the original right-hand side (see \cite{ando2024proceedings} for details on this approach).

For our experiments, we chose $M= 50$, $\lambda = 40$ and $\mu= 2$, and integrated the equation over the time interval $[0,4]$. The dimension of the resulting ODE is thus much larger than the previous examples, yet smaller than many applications which involve the discretization of a PDE. The linear part is dominant but, again contrary to the majority of examples analyzed in the literature of exponential methods, it is not sparse. We compute the error with respect to a reference solution obtained with \texttt{dde23}, a MATLAB integrator for DDEs. The final error, therefore, includes both the contribution of the discretization and that of the numerical method for approximating the solution.

\begin{figure}[ht!]
\centering
\begin{subfigure}{.475\textwidth}
  \centering
  \includegraphics[width=1\linewidth]{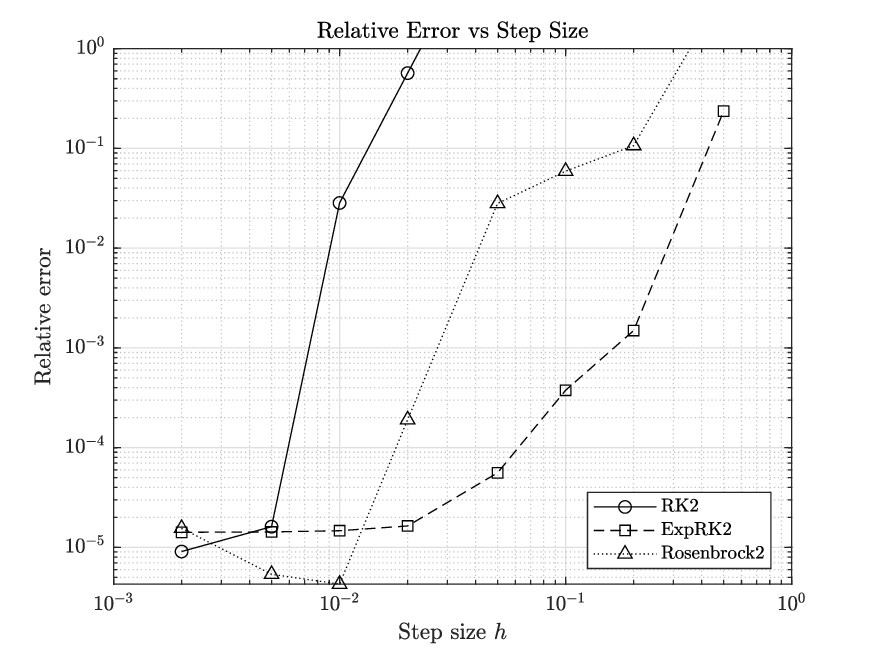}
  \caption{}
  \label{fig:error_step_plotDDE_1_order2}
\end{subfigure}%
\begin{subfigure}{.475\textwidth}
  \centering
  \includegraphics[width=1\linewidth]{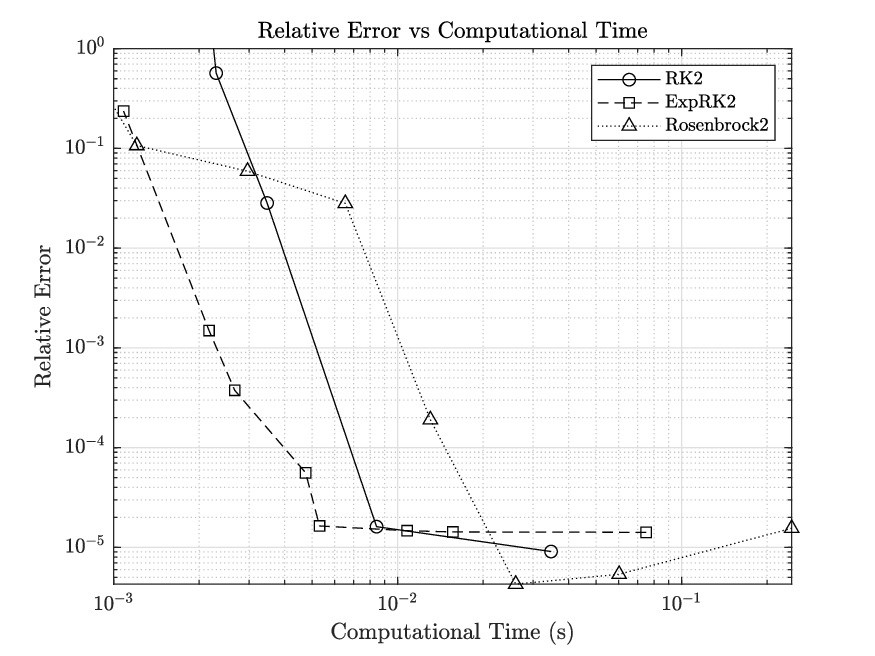}
  \caption{}
  \label{fig:error_time_plotDDE_2_order2}
\end{subfigure}
\begin{subfigure}{.475\textwidth}
  \centering
  \includegraphics[width=1\linewidth]{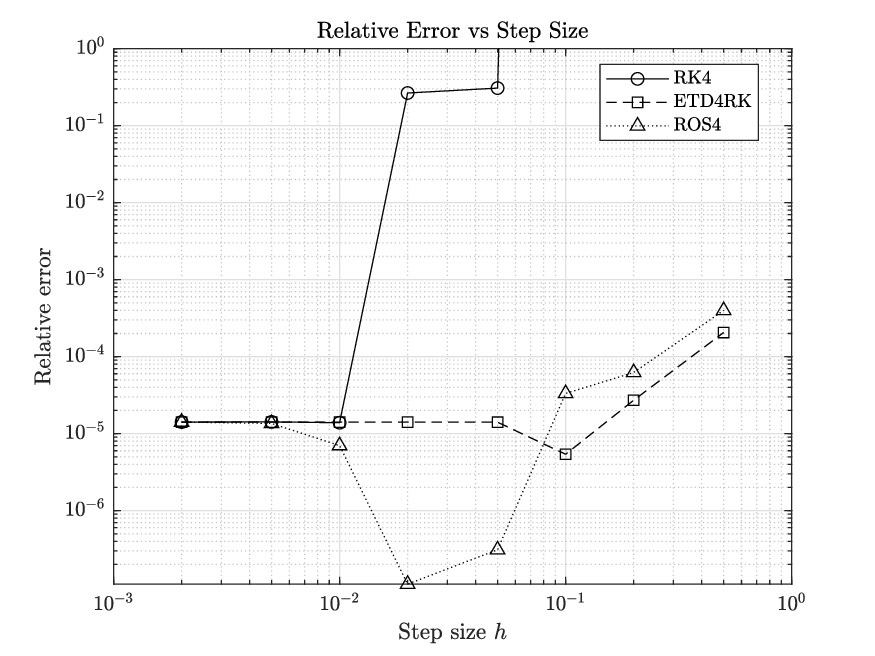}
    \caption{}
  \label{fig:errortime_plotDDE_1_order4}
\end{subfigure}%
\begin{subfigure}{.475\textwidth}
  \centering
  \includegraphics[width=1\linewidth]{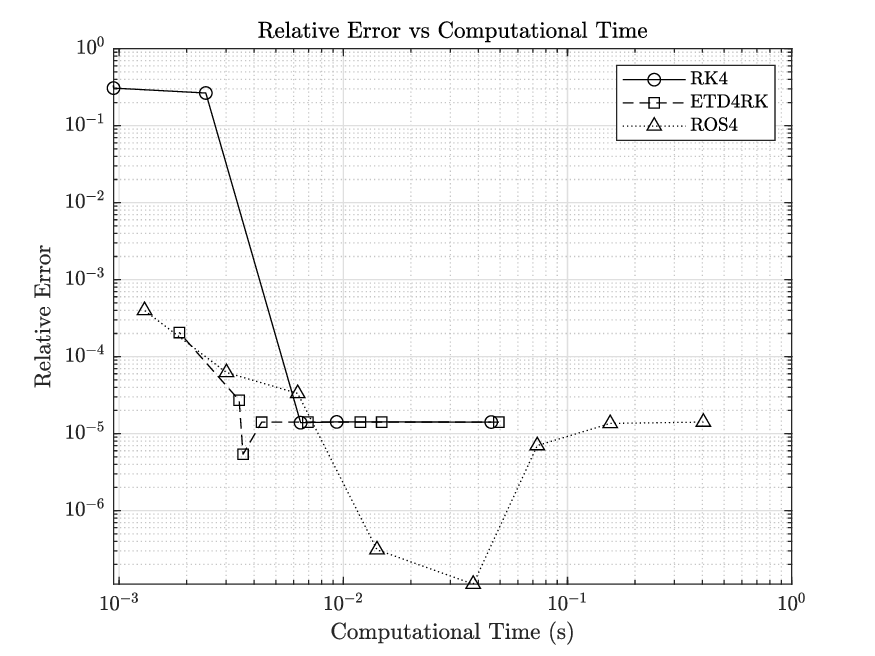}
  \caption{}
\label{fig:errortime_plotDDE_2_order4}
\end{subfigure}
\caption{Relative error as a function of step size, showing how accuracy of each method of order 2
({\subref{fig:error_step_plotDDE_1_order2}}) and order 4 ({\subref{fig:errortime_plotDDE_1_order4}}) evolves with different step choices and relative error as a function of computational time, highlighting the efficiency of methods of order 2 ({\subref{fig:error_time_plotDDE_2_order2}}) and order 4 ({\subref{fig:error_time_plotDDE_2_order2}})}
\label{fig:DDE_plots}
\end{figure}

Consistently with our previous test cases, we present our results in Figure \ref{fig:DDE_plots} and Tables \ref{Table:DDE_order2}-\ref{Table:DDE_o4}. It is immediately clear from all subfigures that the error, for all methods, rarely gets below a certain threshold around $10^{-5}$. We interpret this as the error inherent to the discretization into ODEs and we focus our analysis on the behavior of the methods until they reach the threshold.

The behavior of the classical RK methods is a clear confirmation of the overall stiffness of the problem. Indeed, not only are they affected by stability issues when used with larger step sizes, but their accuracy is consistently lower than the other methods by several order of magnitudes, until all methods reach the aforementioned threshold due to the discretization. The accuracy of the exponential methods is also, in both cases, much higher than the corresponding Rosenbrock methods, although the discrepancy is more evident among the methods of order 2, than of order 4.

In terms of efficiency, the exponential methods are also the ones which reach the threshold first. In the analysis of the order 2 methods, which only reach the threshold with the smallest step sizes, the better efficiency of ExpRK (in terms of required time to reach a certain error) is even clearer.

\paragraph{\bf Final considerations}

The effectiveness of ExpRK2 and ExpRK4 in this example should not come as a surprise, since the chosen value of $\mu$ results in a mildly stiff nonlinear part, whereas the linear term is the main responsible of the overall stiffness. In other words, the problem formulation makes the choice of exponential methods particularly well-suited despite the density of the matrix defining the linear term, which makes the computation of exponential somewhat more onerous.

\begin{table}[ht!]
\centering
\caption{Numerical methods performance for various step sizes (order 2 methods).}
\label{Table:DDE_order2}
\resizebox{\linewidth}{!}{\begin{tabular}{lcccccc}
\toprule
$h$ & RK2 Error & ExpRK2 Error & Rb2 Error & RK2 Time (s) & ExpRK2 Time (s) & Rb2 Time (s) \\
\midrule
$2 \times 10^{-3}$ &$9.0928\times 10^{-6}$ &$	1.4166\times 10^{-5}$ &	$1.5548\times 10^{-5}$ & $0.03471$ &	$0.07505$&$0.24417$\\
$5\times 10^{-3}$&$	1.6142\times 10^{-5}$&$	1.4286\times 10^{-5}$&$	5.3777\times 10^{-6}$&$	0.00844$&$	0.01567$&$	0.06029$\\
$1\times 10^{-2}$&$	2.8393\times 10^{-2}$&$	1.4701\times 10^{-5}$&$	4.2840\times 10^{-6}$&$	0.00347$&$	0.01081$&$	0.02601$\\
$2\times 10^{-2}$&$	5.6934\times 10^{-1}$&$	1.6405\times 10^{-5}$&$	1.9081\times 10^{-4}$&$	0.00230$&$	0.00531$&$	0.01305$\\
$5\times 10^{-2}$&$	2.8932\times 10^{1}$&$	5.5952\times 10^{-5}$&$	2.8142\times 10^{-2}$&$	0.00199$&$	0.00475$&$	0.00655$\\
$1\times 10^{-1}$&$	1.8177\times 10^{29}$&$	3.7608\times 10^{-4}$&$	5.9123\times 10^{-2}$&$	0.00056$&$	0.00267$&$	0.00296$\\
$2\times 10^{-1}$&$	1.8168\times 10^{29}$&$	1.4935\times 10^{-3}$&$	1.0688\times 10^{-1}$&$	0.00028$&$0.00217$&$	0.00121$\\
\bottomrule
\end{tabular}}
\end{table}

\begin{table}[ht!]
\centering
\caption{Numerical methods performance for various step sizes (order 4 methods).}
\label{Table:DDE_o4}
\resizebox{\linewidth}{!}{\begin{tabular}{lcccccc}
\toprule
$h$ & RK4 Error & ExpRK4 Error & Rb4 Error & RK4 Time (s) & ExpRK4 Time (s) & Rb4 Time (s) \\
\midrule
$2 \times 10^{-3}$ & $1.4145 \times 10^{-5}$ & $1.4143 \times 10^{-5}$ & $1.4131 \times 10^{-5}$ & 0.04223 & 0.04528 & 0.34967 \\
$5 \times 10^{-3}$ & $1.4204 \times 10^{-5}$ & $1.4143 \times 10^{-5}$ & $1.3568 \times 10^{-5}$ & 0.01252 & 0.01620 & 0.15586 \\
$1 \times 10^{-2}$ & $1.3871 \times 10^{-5}$ & $1.4143 \times 10^{-5}$ & $6.9969 \times 10^{-6}$ & 0.00397 & 0.00856 & 0.06972 \\
$2 \times 10^{-2}$ & $2.6583 \times 10^{-1}$ & $1.4142 \times 10^{-5}$ & $1.1106 \times 10^{-7}$ & 0.00294 & 0.00741 & 0.03191 \\
$5 \times 10^{-2}$ & $3.0791 \times 10^{-1}$ & $1.4083 \times 10^{-5}$ & $3.1074 \times 10^{-7}$ & 0.00125 & 0.00475 & 0.01593 \\
$1 \times 10^{-1}$ & $1.8129 \times 10^{29}$ & $5.3981 \times 10^{-6}$ & $3.3303 \times 10^{-5}$ & 0.00045 & 0.00412 & 0.00770 \\
$2 \times 10^{-1}$ & $1.4275 \times 10^{42}$ & $2.7167 \times 10^{-5}$ & $6.2284 \times 10^{-5}$ & 0.00024 & 0.00346 & 0.00308 \\
\bottomrule
\end{tabular}}
\end{table}

\section{Conclusions}

In this review, we provided a comprehensive historical and theoretical account of the development of Exponential Runge–Kutta (ExpRK) methods over more than sixty years of research. We highlighted how the methods have evolved from their early formulations in the 1960s to the state-of-the-art developments in stiff order conditions, efficient implementations, and novel applications to partial and delay differential equations. This dual historical and conceptual reconstruction not only situates ExpRK methods within the broader landscape of numerical integration techniques but also clarifies the reasons why they continue to attract substantial attention in both theory and practice.

At the same time, this survey was conceived with a pedagogical aim: to serve as an accessible entry point for readers encountering exponential integrators for the first time. To this end, we have complemented the historical and technical overview with a gradual introduction to the core concepts, supported by carefully selected didactic examples. Our numerical experiments were not intended to showcase the most advanced algorithms or optimized implementations, but rather to illustrate in a transparent and approachable way how ExpRK methods work in practice, and to provide newcomers with concrete tools to grasp their strengths and limitations. In this sense, the manuscript is both a review for specialists and a tutorial for learners. While we focused uniquely on explicit methods for our numerical section, to maintain a coherent presentation, a systematic investigation of their implicit exponential counterparts, together with the study of their robustness properties, remains promising direction for, potentially, a future work.

The dual nature of this contribution, at once a scholarly survey and a didactic introduction, represents, we believe, its distinctive value. On the one hand, specialists may find here a structured synthesis of the major milestones, including recent advances such as high stiff order methods, parallelization strategies, and applications in modeling contexts. On the other hand, students and non-experts may use the text as a stepping stone, moving from the basic principles to more sophisticated topics without losing sight of the underlying intuition. The coexistence of these two layers reflects the current state of the field itself, which remains both historically rich and dynamically evolving.

Looking ahead, it is clear that exponential Runge–Kutta methods will continue to play a central role in computational mathematics. Their capacity to treat stiff linear terms exactly, while retaining the flexibility and structure of Runge–Kutta schemes, ensures their relevance in a wide range of applications, from fundamental scientific modeling to large-scale simulations. Current research trends on stiff order conditions, energy preservation, and efficient implementations for large systems suggest that the methods are still far from reaching a definitive form. On the contrary, they remain a fertile ground for both theoretical inquiry and practical innovation. For these reasons, we expect ExpRK methods to maintain, and even strengthen, their role as a powerful paradigm for the numerical exploration of stiff problems.

%\correz{In this review, we provided a detailed historical report of the evolution in time of the Exponential Runge-Kutta (ExpRK) methods for the numerical integration of (stiff) differential equations. Our focus was on the various ``waves'' of resurfacing research interest in the methods, and the main improvements brought forth by researchers over the years.}

%\correz{In our brief numerical experiments illustrated in a dedicated section, we strove for a fair comparison with existing methods, rather than presenting the best algorithms available on e.g. Matlab. We remark here that, clearly, allowing for a variable step size would likely improve the results we obtain. However, we preferred to present low dimensional systems of differential equations to provide an approachable introduction to the topic to a potential uninitiated reader.}

%\correz{ExpRK methods have been developed and refined for more than sixty years, at the time of writing this survey. Nevertheless, as we have shown, the interest in them has not abated and extensive research on methods of higher stiff order, parallelization and applications to PDEs and DDEs is still ongoing. As they provide a powerful and efficient approach for integrating stiff differential equations, ExpRK methods will likely remain popular for numerical explorations in the near (and not so near) future.}

\textbf{Acknowledgements}\\
\noindent
AA is member and acknowledges the support of {\it Gruppo Nazionale per il Calcolo Scientifico} (GNCS) of {\it Istituto Nazionale di Alta Matematica} (INdAM). Her work was partially supported by the Italian Ministry of University and Research (MUR) through the PRIN 2022 project (No.\ 20229P2HEA) ``Stochastic numerical modelling for sustainable innovation'', Unit of Udine (CUP G53C24000710006). \smallskip

\noindent
NC is member and acknowledges the support of {\it Gruppo Nazionale per l’Analisi Matematica, la Probabilità e le loro Applicazioni} (GNAMPA) of {\it Istituto Nazionale di Alta Matematica} (INdAM) and moreover acknowledges the support of the MIUR - PRIN 2022 project ``Nonlinear dispersive equations in presence of singularities'' (Prot. N. 20225ATSTP).
\smallskip

\noindent MS is member and acknowledges the support of {\it Gruppo Nazionale di Fisica Matematica} (GNFM) of {\it Istituto Nazionale di Alta Matematica} (INdAM). 

\section*{Declarations}
\noindent\textbf{Conflict of interest} The authors declare that they have no conflicts of interest to disclose.

\bibliographystyle{plain}  \small
\bibliography{biblio}% common bib file

\end{document}